\input ./style/arxiv-general.cfg
\documentclass[aos,MSNbibl,nameyear,seceqn,dvips]{arximspdf}
\makeatletter
   \@ifpackageloaded{graphicx}{}{\usepackage{graphicx}}
\makeatother
\usepackage{mathrsfs}
\usepackage{mathbh}

\doi{10.1214/15-AOS1319}
\volume{43}
\issue{4}
\pubyear{2015}
\firstpage{1647}
\lastpage{1681}
\docsubty{FLA}

\makeatletter

\newcommand{\qquand}{ \mbox{and} }

\newcommand{\D}[2][]{\frac{\mathrm{d} #1}{\mathrm{d} #2}}

\newcommand{\indep}{\mathbin{\perp\!\!\!\!\perp}}
\newcommand{\vindep}{\mathbin{ \ddagger}}

\def\mid{|}
\newcommand{\bigmid}{\vert}

\newcommand{\R}{\mathbb{R}}
\newcommand{\Z}{\mathbb{Z}}

\newcommand{\E}{\mathbb{E}}

\newcommand{\cl}{\operatorname{cl}} 
\newcommand{\csep}{\operatorname{sep}} 
\newcommand{\pa}{\operatorname{pa}} 
\newcommand{\an}{\operatorname{an}} 
\newcommand{\pr}{\operatorname{pr}} 

\newcommand{\dist}{\mathscr}
\newcommand{\dnorm}{\dist{N}}
\newcommand{\dinvwish}{\dist{IW}}
\newcommand{\dhinvwish}{\dist{HIW}}

\newcommand{\mobius}{M\"{o}bius }

\newtheorem{thmm}{Theorem}[section]
\newtheorem{lem}[thmm]{Lemma}
\newtheorem{cor}[thmm]{Corollary}
\newtheorem{prp}[thmm]{Proposition}

\newproclaim{dfn}{Definition}[section]
\newproclaim{exm}{Example}[section]
\newproclaim{rem}{Remark}

\newcommand{\erdosrenyi}{Erd\H{o}s--R\'{e}nyi }
\newcommand{\eqref}[1]{(\ref{#1})}
\newcommand{\tmid}{\vert}
\newcommand{\markov}{\mathfrak{M}}
\newcommand{\markoveq}{\stackrel{\markov}{\sim}}
\newcommand{\markovdist}{\mathfrak{P}}
\newcommand{\graph}[1]{\mathcal{#1}}
\newcommand{\G}{\graph{G}}
\newcommand{\rgraph}[1]{\tilde{\graph{#1}}}
\newcommand{\rG}{\rgraph{G}}
\newcommand{\moral}{\mathrm{M}}
\newcommand{\rtheta}{\tilde{\theta}}
\newcommand{\rSigma}{\tilde{\Sigma}}
\newcommand{\rGamma}{\tilde{\Gamma}}
\newcommand{\famtheta}{\vartheta}
\newcommand{\ances}{\mathfrak{D}}
\newcommand{\V}{\mathcal{V}}
\newcommand{\ind}{\mathbh{1}}
\newcommand{\gset}[1]{\mathfrak{#1}}
\newcommand{\F}{\gset{F}}
\newcommand{\sudg}{\gset{U}}
\newcommand{\edge}{\mathcal{E}}
\newcommand{\decom}[3][\sudg]{#1(#2,#3)}
\newcommand{\gtimes}{\otimes}

\newcommand{\dgraph}[1]{\mathcal{#1}}
\newcommand{\dG}{{\dgraph{G}}}
\newcommand{\rdG}{\tilde{\dG}}

\newcommand{\sedag}{\sdag^{\markov}}

\newcommand{\sdag}{\mathfrak{D}}
\newcommand{\sodag}{\mathfrak{D}^{\prec}}
\newcommand{\dagoid}[1]{\mathcal{#1}} 
\newcommand{\sdagoid}{\sdag^{\markov}} 

\renewcommand{\D}{\dagoid{D}}
\newcommand{\rD}{\tilde{\D}}
\newcommand{\law}{\pounds} 
\newcommand{\famlaw}{\mathfrak{L}} 
\newcommand{\glaw}{\mathfrak{G}}
\newcommand{\dprod}{\ltimes}
\newcommand{\dthreegraph}[1]{}

\makeatother

\begin{document}
\begin{frontmatter}

\title{Structural Markov graph laws for Bayesian model~uncertainty}
\runtitle{Structural Markov graph laws}

\begin{aug}
\author[A]{\fnms{Simon} \snm{Byrne}\corref{}\thanksref{T1}\ead[label=e1]{simon.byrne@ucl.ac.uk}}
\and
\author[B]{\fnms{A. Philip} \snm{Dawid}\ead[label=e2]{apd@statslab.cam.ac.uk}}
\runauthor{S. Byrne and A.~P. Dawid}
\affiliation{University College London and University of Cambridge}
\thankstext{T1}{Supported by EPSRC fellowship EP/K005723/1.}

\address[A]{Department of Statistical Science\\
University College London\\
Gower Street, London WC1E 6BT\\
United Kingdom\\
\printead{e1}}

\address[B]{Statistical Laboratory\\
University of Cambridge\\
Wilberforce Road, Cambridge CB3 0WB\\
United Kingdom\\
\printead{e2}}
\end{aug}

%
\received{\smonth{3} \syear{2014}}
%
\revised{\smonth{1} \syear{2015}}

%
\begin{abstract}
This paper considers the problem of defining distributions over graphical
structures. We propose an extension of the hyper Markov properties of
Dawid and Lauritzen [\textit{Ann. Statist.} \textbf{21} (1993) 1272--1317],
which we term \textit{structural Markov properties}, for
both undirected decomposable and directed acyclic graphs, which requires
that the structure of distinct components of the graph be conditionally
independent given the existence of a separating component. This allows the
analysis and comparison of multiple graphical structures, while being able
to take advantage of the common conditional independence
constraints. Moreover, we show that these properties characterise
exponential families, which form conjugate priors under sampling from
compatible Markov distributions.
\end{abstract}

%

%
\end{frontmatter}
%
\section{Introduction}\label{sec1}
A graphical model consists of a graph and a probability distribution
that satisfies a \emph{Markov property} of the graph, being a set of
conditional independence constraints encoded by the graph. Such
models arise naturally in many statistical problems, such as
contingency table analysis and covariance estimation.

\citet{dawid1993} consider distributions over these distributions,
which they term \emph{laws} to emphasise the distinction from the
underlying sampling distribution. Laws arise primarily in two
contexts: as sampling distributions of estimators and as prior and
posterior distributions in Bayesian analyses. Specifically,
\citet{dawid1993} focus on hyper Markov laws that exhibit conditional
independence properties analogous to those of the distributions of the
model. By exploiting such laws, it is possible to perform certain
inferential tasks locally; for instance, posterior laws can be
calculated from subsets of the data pertaining to the parameters of
interest.

Although other types of graphical model exist, we restrict ourselves
to undirected decomposable graphs and directed acyclic graphs, which
exhibit the special property that their Markov distributions can be
constructed in a recursive fashion by taking \emph{Markov
combinations} of smaller components. In the case of undirected
decomposable graphs, for any decomposition $(A,B)$ of the graph $\G$,
a Markov distribution is uniquely determined by the marginal
distributions over $A$ and $B$ [\citet{dawid1993}, Lemma~2.5]. By a
recursion argument, this is equivalent to specifying marginal
distributions on cliques. A similar construction can be derived for
directed acyclic graphs: the distribution of each vertex conditional
on its parent set can be chosen arbitrarily, and the set of such
distributions determines the joint distribution. As we demonstrate in
Section~\ref{sec:markov-equivalence}, this property can also be
characterised in terms of a partitioning based on \emph{ancestral
sets}.

It is this partitioning that makes the notion of hyper Markov laws
possible. In essence, these are laws for which the partitioned
distributions exhibit conditional independence properties analogous to
those of the underlying distributions. In the case of undirected
decomposable graphs, a law $\mbox{\law}$ for $\rtheta$ over
$\markovdist(\G)$, the set of Markov distributions with respect to
$\G$, is \emph{\textup{(}weak\textup{)} hyper Markov} if for any decomposition $(A,B)$,
%
\begin{equation}
\label{eq:whmp-udg} \rtheta_A \indep\rtheta_B \mid
\rtheta_{A \cap B}\qquad [\mbox{\law}].
\end{equation}
Weak hyper Markov laws arise naturally as sampling distributions of
maximum likelihood estimators of graphical models
[\citet{dawid1993}, Theorem~4.22]. A more specific class of laws are those that
satisfy the \emph{strong hyper Markov property}, where for any
decomposition $(A,B)$,
%
\begin{equation}
\label{eq:shmp-udg} \rtheta_{A \mid B} \indep\rtheta_B \qquad [\mbox{\law}].
\end{equation}
When used as prior laws in a Bayesian analysis, strong hyper Markov
laws allow for local posterior updating, in that the posterior law of
clique marginal distributions only depends on the data in the clique
[\citet{dawid1993}, Corollary~5.5].



However, hyper Markov laws only apply to individual graphs: if the
structure of
the graph itself is unknown, then a full Bayesian analysis requires a
prior distribution over graphical structures, which we term a \emph{graph
law}. Very little information is available to guide the choice of such
priors, with a typical choice being a simple uniform or constrained
\erdosrenyi prior.

The aim of this paper is to extend the hyper Markov concept to the structure
of the graph itself. We study graph laws that exhibit similar conditional
independence structure, termed \emph{structural Markov properties}. These
properties exhibit analogous local inference properties, and under minor
assumptions, characterise exponential families, which serve as conjugate
families to families of compatible Markov distributions and hyper
Markov laws.

The outline of the paper is as follows. In
Section~\ref{sec:background} we introduce the terms and notation used
in the paper, in particular the notion of a semi-graphoid to define
what we mean by structure. Section~\ref{sec:udg-smp} develops the
notion of a structural Markov property and characterises such laws for
undirected decomposable graphs. Section~\ref{sec:order-direct-struct}
briefly develops a similar notion for directed graphs consistent with
a fixed ordering. In Section~\ref{sec:markov-equivalence} we consider
the notion of Markov equivalence of directed acyclic graphs, and
extend the structural Markov property to these equivalence classes.
Finally, in Section~\ref{sec:discussion} we discuss some properties,
computational considerations, and future directions.


\section{Background}
\label{sec:background}

Much of the terminology in this paper is standard in the graphical
modelling literature. For this we refer the reader to texts such as
\citet{lauritzen1996} or \citet{cowell2007}. For clarity and
consistency, the following presents some specific terms and notation
used in this paper.

\subsection{Graphs}
\label{sec:graphs}

A \emph{graph} $\G$ consists of a set of \emph{vertices} $\V(\G)$ and
a set of \emph{edges} $\edge(\G)$ of pairs of vertices. In the case
of undirected graphs, $\edge(\G)$ will be a set of unordered pairs of
vertices $\{u,v\}$; in the directed case it will be a set of ordered
pairs $(u,v)$, denoting an arrow from $u$ to $v$, of which $v$ is
termed the \emph{head}. For any subset $A
\subseteq\V(\G)$, $\G_A$ will denote the induced subgraph with vertex
set $A$. A graph is \emph{complete} if there exists an edge between
every pair of vertices, and \emph{sparse} if no edges are present (the
graph with empty vertex set is both complete and sparse).

We focus on two particular classes of graphs.

\subsubsection{Undirected decomposable graphs}
\label{sec:undir-decomp-graphs}

A \emph{path} in an undirected graph $\G$ is a sequence of vertices
$v_0,v_1,\ldots,v_k$ such that $\{v_i,v_{i+1}\} \in\edge(\G)$, in
which case we can say $v_0$ is \emph{connected} to $v_k$. Sets $A,B
\subseteq\V(\G)$ are \emph{separated} by $S \subseteq\V(\G)$ if
every path starting at an element of $A$ and ending at an element of
$B$ contains an element of $S$.

A pair of sets $(A,B)$ is a \emph{covering pair} of $\G$ if $A \cup B =
\V(\G)$. A covering pair is a \emph{decomposition} if $\G_{A \cap
B}$ is
complete and $A$ and $B$ are separated by $A \cap B$ in $\G$. A decomposition
is \emph{proper} if both $A$ and $B$ are strict subsets of $\V(\G)$.
For any set
of undirected graphs $\F$, define $\decom[\F]{A}{B}$ to be the set
of $\G\in
\F$ for which $(A,B)$ is a decomposition.

A graph is \emph{decomposable} if it can be recursively decomposed
into complete subgraphs. An equivalent condition is that the graph is
\emph{chordal}, in that there exists no set which induces a cycle
graph of length 4 or greater. Throughout the paper we will take $V$
to be a fixed, finite set, and define $\sudg$ to be the set of
undirected decomposable graphs (UDGs) with vertex set $V$.

The maximal sets inducing complete subgraphs are termed
\emph{cliques}, the set of which is denoted by $\cl(\G)$. For any
decomposable graph it is possible to construct a \emph{junction tree}
of the cliques. The intersections of neighbouring cliques in a
junction tree are termed \emph{\textup{(}clique\textup{)} separators}, the set of which
is denoted by $\csep(\G)$. The \emph{multiplicity} of a separator is
the number of times it appears in the junction tree. The cliques,
separators, and their multiplicities are invariants of the graph.

An undirected graph $\G$ is \emph{collapsible} onto $A \subseteq
\V(\G)$ if each connected component $C$ of $\G_{V \setminus A}$ has a
boundary $B = \{u \dvtx \{u,v\} \in\edge(\G), v \in C, u \notin C\}$
which induces a complete subgraph. Note that if $(A,B)$ is a
decomposition of $\G$, then $\G$ is collapsible onto both $A$ and $B$.

\subsubsection{Directed acyclic graphs}
\label{sec:direct-acycl-graphs}

A directed graph $\G$ is \emph{acyclic} if there exists a
\emph{compatible well-ordering} $\prec$ on $\V(\G)$, that is, such
that $u
\prec v$ for all $(u,v) \in\edge(\G)$. For any such $\prec$, the
\emph{predecessors} of a vertex $v$ is the set $\pr_\prec(v) = \{u
\in
\V(\G) \dvtx u \prec v)$. The set of directed acyclic graphs (DAGs) on
$V$ will be denoted by $\sdag$, and the subset for which $\prec$ is a
compatible well-ordering is denoted by $\sodag$.

A vertex $u$ is a \emph{parent} of $v$ if $(u,v) \in\edge(\G)$. The
set of parents of $v$ is denoted by $\pa_\G(v)$. Conversely, $u$ is a
\emph{child} of $v$. A set $A \subseteq\V(\G)$ is \emph{ancestral}
in $\G$ if $v \in A\Rightarrow\pa_\G(v) \subseteq A$. The minimal
ancestral set containing $B \subseteq\V(\G)$ is denoted by
$\an_\G(B)$.

The \emph{skeleton} of a directed graph $\G$ is the undirected graph
obtained by replacing all the directed edges with undirected edges.
The \emph{moral graph} of $\G$, denoted by $\dG^{\moral}$, is the
skeleton of the graph obtained by adding (if necessary) an edge
between each pair of vertices having a common child.

\subsection{Distributions and laws}
\label{sec:graphical-models}

Let $X = (X_v)_{v \in V}$ be a random vector on some product space
$\prod_{v \in V} \mathcal{X}_v$, with distribution denoted by $P$ or
$\theta$. A \emph{model} is a family of distributions $\Theta$ for
$X$.

Following \citet{dawid1993}, a distribution over $\Theta$ will be
termed a \emph{law} and denoted by $\mbox{\law}$. A random distribution
following such a law will be denoted by $\rtheta$.

For any $A \subseteq V$, $X_A$ will denote the subvector $(X_v)_{v \in
V}$, with $P_A$ or $\theta_A$ denoting its marginal distribution.
The marginal law of $\rtheta_A$ will be denoted by $\mbox{\law}_A$.
Furthermore, for any pair $A,B \subseteq V$, we can denote by
$\theta_{A \mid B}$ the collection of conditional distributions of
$X_A \mid X_B$ under $\theta$, and by $\mbox{\law}_{A \mid B}$ the induced
law of $\rtheta_{A \mid B}$ under $\mbox{\law}$. We will use $\simeq$ to
indicate the existence of a bijective function; for instance, we can
write $(\theta_A, \theta_{V \mid A}) \simeq\theta$ for any $A
\subseteq V$.

\subsection{Semi-graphoids}
\label{sec:separoids}

When discussing the ``structure'' of a graphical\break model, many authors
use this term to refer to the graph itself. In particular, when they
talk of ``estimating the structure,'' they mean inferring the presence
or absence of individual edges of the graph.

In this paper, we take the view that ``structure'' refers to a set of
conditional independence properties, and that a graph is merely a
representation of this structure. This distinction is an important
one: it implies that graphs that encode the same set of conditional
independence statements must be treated as identical, leading to the
notion of \emph{Markov equivalence}. A more subtle but even more
important point is that when investigating properties such as
decompositions or ancestral sets, we are, effectively, looking at
properties of sets of conditional independencies.

To make this more concrete, we use the notion of a
\emph{semi-graphoid}, a special case of a \emph{separoid}
[\citet{dawid2001}], to describe the abstract properties of conditional
independence.
%
\begin{dfn}
\label{dfn:separoid}
Given a finite set $V$, a \emph{semi-graphoid} is a set $M$ of
triples of the form $\langle A, B \mid C \rangle$, where $A,B,C
\subseteq V$, satisfying the properties:
\begin{longlist}[S0]
\item[S0] for all $A,B \subseteq V$, $\langle A, B \mid A \rangle
\in M$;
\item[S1] if $\langle A,B \tmid C \rangle\in M$, then $\langle B,A
\tmid C \rangle\in M$;
\item[S2] if $\langle A,B \tmid C \rangle\in M$ and $D \subseteq
A$, then $\langle D,B \tmid C \rangle\in M$;
\item[S3] if $\langle A,B \tmid C \rangle\in M$ and $D \subseteq
A$, then $\langle A,B \tmid C \cup D \rangle\in M$;
\item[S4] if $\langle A,B \tmid C \rangle\in M$ and $\langle A,D
\tmid B \cup C \rangle\in M$, then $\langle A,B \cup D \tmid C
\rangle\in M$.
\end{longlist}
\end{dfn}
These properties match the well-established properties of conditional
independence [\citet{dawid1979}].

We can define the semi-graphoid of a graph as the set of triples
encoding its global Markov property: the semi-graphoid of an
undirected graph $\G$ is
%
\begin{equation}
\label{eq:undirected-global} \markov(\G) = \bigl\{ \langle A,B \tmid C \rangle \dvtx \mbox{$A$
and $B$ are separated by $C$ in $\G$} \bigr\},
\end{equation}
and the semi-graphoid of a directed acyclic graph $\dG$ is the set
%
\begin{equation}
\label{eq:directed-global} \markov(\dG) = \bigl\{ \langle A,B \tmid C \rangle \dvtx
\mbox{$A$ and $B$ are separated by $C$ in $\dG^\moral_{\an(A \cup
B \cup C)}$}
\bigr\}.
\end{equation}

We say that a joint distribution $P$ for $X = (X_v)_{v \in V}$ is
\emph{Markov} with respect to a semi-graphoid $M$ if
\[
\langle A, B \tmid C \rangle\in M \quad\Rightarrow\quad X_A \indep
X_B \mid X_C\qquad  [P].
\]
That is, a distribution is Markov with respect to a graph if it is
Markov with respect to the semi-graphoid of the graph. We write
$\markovdist(\G)$ or $\markovdist(M)$ to be the set of distributions
that are Markov with respect to $\G$ or $M$.

Similarly, a law $\mbox{\law}$ is \emph{weak hyper Markov} with respect to
the semi-graphoid if
\[
\langle A, B \tmid C \rangle\in M \quad\Rightarrow\quad \rtheta_{A \cup C} \indep
\rtheta_{B \cup C} \mid\rtheta_C\qquad [\mbox{\law}].
\]
However, the strong hyper Markov laws cannot be directly characterised
in terms of the semi graphoid.

Semi-graphoids have a natural projection operation: for any set $U
\subseteq V$, we can define the projection onto $U$ of a semi-graphoid
$M$ on $V$ to be
\[
M_U = \bigl\{ \langle A, B \mid C \rangle\in M \dvtx A,B,C \subseteq
U \bigr\}.
\]
Under certain conditions, this can match the natural projection
operation, the induced subgraph, of the underlying graph. For
undirected graphs,\break $[\markov(\G)]_U = \markov(G_U)$ if and only if
$\G$ is collapsible onto $U$ [\citet{asmussen1983}, Corollary~2.5]. For
directed acyclic graphs, we have the weaker sufficient condition that
if $A$ is ancestral in $\dG$, then $[\markov(\dG)]_A =
\markov(\dG_A)$.


\section{Undirected structural Markov property}
\label{sec:udg-smp}

We now extend the hyper Markov framework to the case where the graph
itself is regarded as a random object $\rG$, taking values in the set
of undirected decomposable graphs with vertex set $V$; equivalently,
$\rG$ can be thought of as a random vector of length ${|V|\choose2}$
indicating the presence or absence of individual edges. As the graph
is a parameter of the model, we term its distribution a \emph{graph
law}, denoted by $\glaw(\rG)$. Our aim is to identify and
characterise hyper Markov-type properties for $\rG$.

Hyper Markov laws are motivated by the property that graph
decompositions allow one to decompose Markov distributions into
separate components. For a fixed graph $\G\in\decom{A}{B}$, then
any Markov distribution $\theta\in\markovdist(\G)$ is uniquely
characterised by its marginals $\theta_A$ and $\theta_B$, taking
values in $\markovdist(\G_A)$ and $\markovdist(\G_B)$, respectively
[\citet{dawid1993}, Lemma~2.5]. Moreover, these can be chosen
arbitrarily, subject only to the constraint $(\theta_A)_{A \cap B} =
(\theta_B)_{A \cap B}$. Hyper Markov laws are derived by imposing
probabilistic conditional independence on this natural separation.

In a similar manner, graphs themselves can be characterised by their
projections onto each part of a decomposition.
%
\begin{prp}
\label{prp:graph-product}
Let $\graph{H}$ and $\graph{J}$ be decomposable graphs with vertex
set $A$ and $B$, respectively, such that both $\graph{H}_{A \cap B}$
and $\graph{J}_{A \cap B}$ are complete. Then the graph $\G$ with
$\edge(\G) = \edge(\graph{H}) \cup\edge(\graph{J})$ is the unique
decomposable graph on $A \cup B$ such that:
\begin{longlist}[(iii)]
\item[(i)] $\G_A = \graph{H}$,
\item[(ii)] $\G_B = \graph{J}$, and
\item[(iii)] $(A,B)$ is a decomposition of $\G$.
\end{longlist}
\end{prp}
\begin{pf}
To satisfy (i) and (ii), the edge set must contain $\edge(\graph{H})
\cup\edge(\graph{J})$. It cannot contain any additional edges
$\{u,v\}$, as this would violate: (i), if $\{u,v\} \subseteq A$;
(ii), if $\{u,v\} \subseteq B$; or (iii), if $u \in A \setminus B$
and $v \in B \setminus A$.
\end{pf}
In other words, a graph $\G\in\decom{A}{B}$ is characterised by
$\G_A$ and $\G_B$, and $\G_A$ and $\G_B$ can be chosen independently.
Moreover, this also decomposes the semi-graphoid, as $\G$ is
collapsible onto both $A$ and $B$.

We define the graph $ G$ resulting from
Proposition~\ref{prp:graph-product} to be the \emph{graph product} of
$\graph{H}$ and $\graph{J}$, denoted by
\[
\G= \graph{H} \gtimes\graph{J}.
\]
%
\begin{rem*}
Although we only use the graph product when $\G_{A \cap B}$ is
complete, the
definition can be extended to the case where $\graph{H}$ and $\graph
{J}$ are
collapsible onto $A \cap B$.
\end{rem*}

For a graph law $\glaw(\rG)$ over $\decom{A}{B}$, a straightforward
way to extend the hyper Markov property in this case would be to
require that
%
\begin{equation}
\label{eq:smp-sub} \rG_A \indep\rG_B \mid
\rG_{A \cap B}\qquad [\glaw].
\end{equation}
Note that in this case the term $\rG_{A \cap B}$ is redundant: if
$(A,B)$ is a decomposition of~$\G$, then $\G_{A \cap B}$ must be
complete, and so we are left with a statement of marginal independence
$\rG_A \indep\rG_B$.

A more general question remains: how might this property be extended
to a graph law over all undirected graphs? A seemingly simple
requirement is that \eqref{eq:smp-sub} should hold whenever a
decomposition exists. This motivates the following definition.
%
\begin{dfn}[(Structural Markov property)]
\label{dfn:smp}
A graph law $\glaw(\rG)$ over $\sudg$ is \emph{structurally Markov}
if for any covering pair $(A,B)$ where $\glaw(\decom{A}{B}) > 0$,
then $\rG_A$ is independent of $\rG_B$, conditional on $(A,B)$ being
a decomposition of $\rG$. This is written as
%
\begin{equation}
\label{eq:smp} \rG_A \indep\rG_B \bigmid\bigl\{\rG\in
\decom{A} {B}\bigr\}\qquad [\glaw].
\end{equation}
\end{dfn}
%

%
%
%
%
%
%
%
%

\begin{figure}

\includegraphics{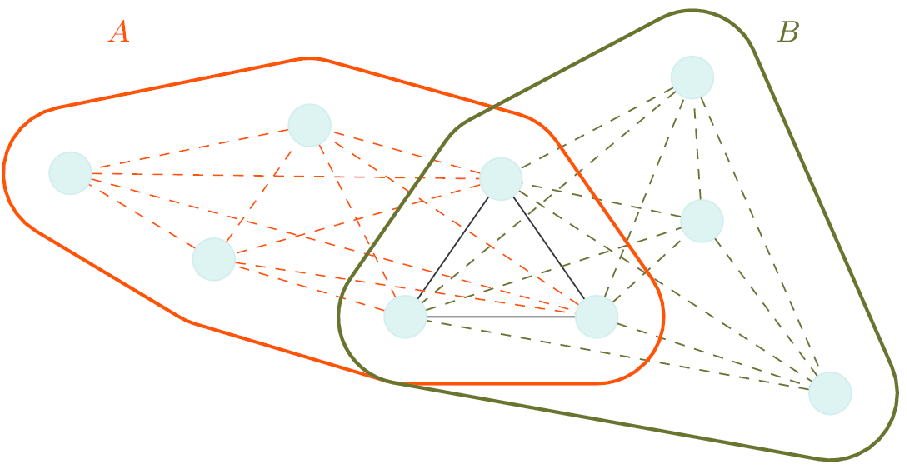}

\caption{A representation of the structural Markov
property for undirected graphs. Conditional on $(A,B)$ being a
decomposition, the existence of the remaining edges in $\rG_A$
(\protect
\includegraphics{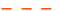}) is independent of those in
$\rG_B$ (\protect
\includegraphics{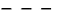}).}
\label{fig:udg-smp}
\end{figure}

In essence, the structural Markov property states that the structures
of different induced subgraphs are conditionally independent given
that they are in separate parts of a decomposition. See Figure~\ref{fig:udg-smp} for a depiction.

The use of braces on the right-hand side of \eqref{eq:smp} is to
emphasise that the conditional independence is defined with respect to
the \emph{event} $\rG\in\decom{A}{B}$, and not a random variable as
in the Markov and hyper Markov properties. In other words, we do not
assume $\rG_A \indep\rG_B \mid\{\rG\notin\decom{A}{B}\}$.

\subsection{Products and projections}
\label{sec:graph-product}

The graph product operation provides a very useful characterisation of
the structural Markov property.
%
\begin{prp}
\label{prp:sm-ratio}
A graph law $\glaw$ is structurally Markov if and only if for every
covering pair $(A,B)$, and every $\G,\G' \in\decom{A}{B}$,
%
\begin{equation}
\label{eq:sm-ratio} \pi(\G) \pi\bigl(\G'\bigr) = \pi\bigl(
\G_A \gtimes\G'_B\bigr) \pi\bigl(
\G'_A \gtimes\G_B\bigr),
\end{equation}
where $\pi$ is the density of $\glaw$ with respect to the counting
measure on $\sudg$.
\end{prp}
\begin{pf}
By Proposition~\ref{prp:graph-product}, both $\G_A \gtimes\G'_B, \G
'_A \gtimes\G_B \in
\decom{A}{B}$, and so if $\glaw(\decom{A}{B}) = 0$, the statement is trivial.
Otherwise, the conditional density of a structural Markov law is of
the form
\[
\pi\bigl(\G\mid\bigl\{\G\in\decom{A} {B}\bigr\}\bigr) = \pi\bigl(\G_A \mid
\decom{A} {B}\bigr) \pi\bigl(\G_B \mid\decom{A} {B}\bigr).
\]
The result follows by substitution into \eqref{eq:sm-ratio}.
\end{pf}

The structural Markov property has an inherent divisibility property
that arises on subgraphs induced by decompositions. First we require
the following lemma.
%
\begin{lem}
\label{lem:sub-decomposition}
Let $(A,B)$ be a decomposition of a graph $\G$, and $(S,T)$ a
covering pair of $A$ with $A \cap B \subseteq T$. Then $(S,T)$ is a
decomposition of $\G_A$ if and only if $(S,T \cup B)$ is a
decomposition of $\G$.
\end{lem}
\begin{pf}
Recall that $W$ separates $U$ and $V$ in $\G$ if and only if
$\langle U,V \tmid W \rangle\in\markov(\G)$. Since $(S,T)$ is a
covering pair of $A$, $\langle S \cup T, B \tmid S \cap B \rangle
\in\markov(\G)$, and hence $\langle S, B \tmid T \rangle\in
\markov(\G)$. If $(S,T)$ is a decomposition of $G_A$, then $\langle
S, T \tmid S \cap T \rangle\in\markov(\G_A)$, which implies that
$\langle S, B \cup T \tmid T \cap S \rangle\in\markov(\G)$. Since
$\G_{(S \cup B) \cap T} = \G_{T \cap S}$ is complete, $(S \cup B,
T)$ is a decomposition of $\G$.

The converse result follows by the reverse argument.
\end{pf}
%
\begin{thmm}
\label{thmm:sm-subgraph}
Let $\glaw(\rG)$ be a structurally Markov graph law. Then the
conditional law for $\rG_A \mid\{\rG\in\decom{A}{B}\}$ is also
structurally Markov.
\end{thmm}
\begin{pf}
Let $(S,T)$ be a covering pair of $A$: If we restrict $\rG\in
\decom{A}{B}$, then $\rG_{A \cap B}$ must be complete. As we are
only interested in the case where $(S,T)$ is a decomposition of
$\rG_A$, then $A \cap B$ must be a subset of either $S$ or $T$:
without loss of generality, we may assume $A \cap B \subseteq T$.

Since $(S,T \cup B)$ is a covering pair of $V$, by the structural
Markov property,
\[
\rG_S \indep\rG_{T \cup B} \bigmid\bigl\{\rG\in\decom{S} {T \cup B}
\bigr\}.
\]
If $\ind_{E}$ is the indicator variable of an event $E$, we can
write
\[
\rG_S \indep(\rG_T, \ind_{\rG_{T \cup B} \in\decom{T}{B}}) \bigmid\bigl\{\rG
\in\decom{S} {T \cup B}\bigr\}.
\]
By the properties of conditional independence [\citet{dawid1979}], the
term $\ind_{\rG_{T \cup B} \in\decom{T}{B}}$ may be moved to the
right-hand side. Furthermore, we are only interested in the case
where it equals $1$. Hence we can write
\[
\rG_S \indep\rG_T \bigmid\bigl\{\G_{T \cup B} \in
\decom{T} {B}\bigr\}, \bigl\{ \rG\in\decom{S} {T \cup B}\bigr\}.
\]
By Lemma~\ref{lem:sub-decomposition}, $\rG_{T \cup B} \in
\decom{T}{B}$ if and only if $\rG\in\decom{S \cup T}{B} =
\decom{A}{B}$. So
\[
\rG_S \indep\rG_T \bigmid\bigl\{\rG\in\decom{A} {B}\bigr\}, \bigl\{
\rG\in \decom{S} {T \cup B}\bigr\}.
\]
Again, by Lemma~\ref{lem:sub-decomposition}, $\rG\in\decom{S}{T
\cup B}$ if and only if $\rG_A \in\decom{S}{T}$, hence:
\[
\rG_S \indep\rG_T \bigmid\bigl\{\rG\in\decom{A} {B}\bigr\},
\bigl\{
\rG_A \in \decom{S} {T}\bigr\}. 
\]
\upqed\end{pf}

\subsection{Structural meta Markov property}
\label{sec:meta-struct-mark}

Dawid and Lauritzen (\citeyear{dawid1993})
define a \emph{meta Markov model} as a set of Markov
distributions that exhibits \emph{conditional variation independence},
denoted by the ternary relation ($\cdot\vindep\cdot\mid\cdot$), in
place of the conditional probabilistic independence of hyper Markov
laws; see also \citet{dawid2001b}. Analogous structural properties
can be defined for families of graphs.

\begin{dfn}[(Structural meta Markov property)]
\label{dfn:struct-meta-markov}
Let $\F$ be a family of undirected decomposable graphs on $V$. Then
$\F$ is \emph{structurally meta Markov} if for every covering pair
$(A,B)$, the set $\{\G_A \dvtx \G\in\decom[\F]{A}{B}, \G_B =
\graph{J}\}$ is the same for all $\graph{J} \in
[\decom[\F]{A}{B}]_B]$. That is,
\[
\G_A \vindep\G_B \bigmid\bigl\{\G\in\decom[\F]{A} {B}
\bigr\}.
\]
\end{dfn}
In other words, this property requires that the set of pairs
$(\G_A,\G_B)$ of $\G\in\decom[F]{A}{B}$ be a product set. Clearly
the set $\sudg$ of all decomposable graphs on $V$ is structurally meta
Markov.

As with probabilistic independence, we can characterise it in terms of
the graph product operation.
%
\begin{thmm}
\label{thmm:struct-meta-markov-gtimes}
A family of undirected decomposable graphs $\F$ is structurally meta
Markov if and only if $\G_A \gtimes\G'_B \in\F$ for all $\G,\G'
\in\decom[\F]{A}{B}$.
\end{thmm}
\begin{pf}
This follows directly from Proposition~\ref{prp:graph-product}.
\end{pf}
Theorem~\ref{thmm:struct-meta-markov-gtimes} is particularly useful in
that if a family of graphs is characterised by a specific property, we
can show that it is structurally meta Markov if this property is
preserved under the graph product operation.

\begin{exm}
\label{exm:struct-meta-markov-maximal-clique}
The set of undirected decomposable graphs whose clique size is
bounded above by some $n$ and whose separator size is bounded below
by $m$ is structurally meta Markov. To see this, note that a clique
of $\G_A \gtimes\G'_B$ must be a clique of either $\G_A$ or $\G'_B$
(and hence of either $\G$ or $\G'$), and therefore the graph product
operation cannot increase the size of the largest clique.
Similarly, it is not possible for a graph product to decrease the size of
the smallest separator: a separator of $\G_A \gtimes\G'_B$ must
either be a
separator of $\G$ or $\G'$, or be $A \cap B$ (this is a consequence of
Lemma~\ref{lem:ctdecom}).

In the case $n=2$ and $m=0$, this is the set of forests on $V$, and
when $n=2$ and $m=1$, this is the set of trees on $V$.
\end{exm}
%
\begin{exm}
\label{exm:struct-meta-markov-sandwich}
Consider two graphs $\G^L,\G^U \in\sudg$ such that $\edge(\G^L)
\subseteq\edge(\G^U)$. Then the ``sandwich'' set between the two
graphs,
\[
\bigl\{ \G\in\sudg\dvtx \edge\bigl(\G^L\bigr) \subseteq\edge(\G)
\subseteq\edge \bigl(\G^U\bigr) \bigr\},
\]
is structurally meta Markov. This follows from the fact that an
edge can only appear in a graph product if it is in one of the
elements of the product.
\end{exm}

As with hyper Markov laws, being a structural meta Markov family is a
necessary condition for the existence of a structural Markov law.
%
\begin{thmm}
\label{thmm:struct-markov-support-meta}
The support of a structurally Markov graph law is a structurally
meta Markov family.
\end{thmm}
\begin{pf}
Let $\F$ be the support of the structurally Markov graph law $\glaw$
with density $\pi$. By Proposition~\ref{prp:sm-ratio}, if $\G,\G'
\in\decom[\F]{A}{B}$ and both $\pi(\G)$ and $\pi(\G')$ are
nonzero, then $\pi(\G_A \gtimes\G'_B)$ must also be nonzero, and
hence in $\decom[\F]{A}{B}$. Therefore, by Theorem~\ref{thmm:struct-meta-markov-gtimes}, $\F$ is structurally meta
Markov.
\end{pf}

\subsection{Compatible distributions and laws}

We now investigate how the structural Markov property interacts with
the Markov and hyper Markov properties. In order to do this, we need
to define families of distributions and laws for every graph.

\begin{dfn}
\label{dfn:compat}
For $\F\subseteq\sudg$, let $\famtheta=\{\theta^{(\G)}\dvtx \G\in\F
\}$ be a family of probability distributions for $X$. We write $X
\sim\famtheta\mid\rG$ if, given $\rG=\G$, $X \sim\theta^{(\G)}$.
Then $\famtheta$ is \emph{compatible} if:
\begin{longlist}[(ii)]
\item[(i)] for each $\G\in\F$, $X$ is Markov with respect to $\G$
under $\theta^{(\G)}$, and
\item[(ii)] $\theta_A^{(\G)} = \theta_A^{(\G')}$ whenever $\G,\G'
\in\F$ are collapsible onto $A$ and $\G_A = \G'_A$.
\end{longlist}
\end{dfn}
Similar properties can be defined for laws.
%
\begin{dfn}
\label{dfn:hyper-compat}
For $\F\subseteq\sudg$, let $\famlaw= \{\mbox{\law}^{(\G)} \dvtx \G\in
\F\}$ be a family of laws for the parameters $\rtheta$ of a family
of distributions on $X$. Again, we can write $\rtheta\sim\famlaw
\mid\rG$ if, given $\rG=\G$, $\rtheta\sim\mbox{\law}^{(\G)}$. Then
$\famlaw$ is \emph{hyper compatible} if:
\begin{longlist}[(ii)]
\item[(i)] for all $\G\in\F$, $\mbox{\law}^{(\G)}$ is weak hyper Markov
with respect to $\G$, and
\item[(ii)] $\mbox{\law}^{(\G)}_A = \mbox{\law}^{(\G')}_A$ whenever $\G,\G'\in
\F$
are collapsible onto $A$ and $\G_A = \G'_A$.
\end{longlist}
\end{dfn}

\begin{rem*}
\citet{dawid1993}, Section~6.2, originally used the term compatible
to refer to what we term the hyper compatible case: we introduce the
distinction so as to extend the terminology to the distributional
(nonhyper) case.
\end{rem*}

As Markov distributions and hyper Markov laws are characterised by
their clique-marginal distributions
[\citet{dawid1993}, Theorems 2.6 and
3.9], it is sufficient for condition (ii) in Definitions
\ref{dfn:compat} and \ref{dfn:hyper-compat} to hold when $\G_A$ and
$\G'_A$ are complete. Moreover, if the complete graph $\G^{(V)}$ is
contained in $\F$, then the compatible and hyper compatible families
are characterised entirely by $\theta^{(\G^{(V)})}$ and
$\mbox{\law}^{(\G^{(V)})}$, respectively.

\begin{exm}
\label{exm:hyper-inv-wishart}
The \emph{inverse Wishart law} for the covariance selection model
$\theta(X) = \dnorm(0,\Sigma)$ assigns $\mbox{\law}(\Sigma) =
\dinvwish(\delta; \Phi)$. This law is strong hyper Markov with
respect to the complete graph on $V$, and the hyper compatible
family generated by $\mbox{\law}$ are the \emph{hyper inverse Wishart laws}
$\mbox{\law}^{(\G)}(\Sigma) = \dhinvwish_{\G} (\delta; \Phi)$
[\citet{dawid1993}, Example~7.3].
\end{exm}

A law induces marginal distribution $\theta_{\mbox{\scriptsize{\law}}}$ for $X$ such that
$\theta_{\mbox{\scriptsize{\law}}}(A) = \E_{\mbox{\scriptsize{\law}}}[\rtheta(A)]$, referred to as the
\emph{predictive distribution} in Bayesian problems. Therefore a
family of laws will also induce a family of distributions. Although
in general hyper compatibility will not imply compatibility, there is
one important special case.
%
\begin{prp}
\label{prp:compat-law-dist}
Let $\famlaw$ be a family of laws such that each law $\mbox{\law}^{(\G)}
\in\famlaw$ is strong hyper Markov. Then the family of marginal
distributions
\[
\{ \theta_{\mbox{\scriptsize{\law}}}\dvtx \mbox{\law}\in\famlaw\}
\]
is hyper compatible.
\end{prp}
\begin{pf}
By \citet{dawid1993}, Proposition~5.6, the marginal distribution of
a strong hyper Markov law is Markov with respect to the same graph.
The result follows by noting that the marginal distribution on a
complete subgraph is a function of the marginal law.
\end{pf}

A graph law $\glaw(\rG)$ combined with a compatible set of
distributions $\famtheta$ defines a joint distribution $(\glaw,
\famtheta)$ for $(\rG,X)$ under which $X \mid\rG=\G\sim
\theta^{(\G)}$. Likewise, $\glaw$ combined with a set of hyper
compatible laws $\famlaw$ defines a joint law $(\glaw,\famlaw)$ for
$(\rG, \rtheta)$, and so a joint distribution on $(\rG, \rtheta, X)$.

The key conditional independence property of any such joint
distribution or law can be characterised as follows.
%
\begin{prp}
\label{prp:compat-ci}
For any graph law $\glaw$ over $\F\subseteq\sudg$ for $\rG$, and
$X \sim\famtheta$ for a compatible family $\famtheta$ indexed by
$\F$,
\[
X_A \indep\rG_B \bigmid\rG_A, \bigl\{\rG\in
\decom{A} {B}\bigr\}\qquad [\glaw, \famtheta].
\]
Similarly, if $\rtheta\sim\famlaw$ for a hyper compatible family
$\famlaw$ indexed by $\F$, then
\[
\rtheta_A \indep\rG_B \bigmid\rG_A, \bigl\{\rG
\in\decom{A} {B}\bigr\}\qquad [\glaw, \famlaw].
\]
\end{prp}
\begin{pf}
Let $\G,\G' \in\decom{A}{B}$ such that $\G_A = \G'_A$. As $\G$ and
$\G'$ are both collapsible onto $A$, then $\theta^{(\G)}_A =
\theta^{(\G')}_A$ in a compatible family, and $\mbox{\law}^{(\G)}_A =
\mbox{\law}^{(\G')}_A$ in a hyper compatible family.\vadjust{\goodbreak}
\end{pf}

When combined with the structural Markov property, we obtain some
useful results.
%
\begin{thmm}
\label{thmm:sm-mark}
If $\rG$ has a structurally Markov graph law $\glaw$, and $X$ has a
distribution from a compatible set $\famtheta$, then
\[
(X_A, \rG_A) \indep(X_B,
\rG_B) \bigmid X_{A \cap B}, \bigl\{\rG\in\decom{A} {B}\bigr\}\qquad [\glaw,
\famtheta].
\]
\end{thmm}
\begin{pf}
See Appendix \ref{sec:proofs}.
\end{pf}

\begin{cor}
\label{cor:sm-mark-post}
If $\rG$ has a structurally Markov graph law, and $X$ has a
distribution from a compatible set $\famtheta$, then the posterior
graph law for $\rG$ is structurally Markov.
\end{cor}
\begin{pf}
By Theorem~\ref{thmm:sm-mark} and the axioms of conditional
independence, we easily obtain
\[
\rG_A \indep\rG_B \bigmid X, \bigl\{\rG\in\decom{A} {B}\bigr\}.
\]
\upqed\end{pf}

We can also apply similar arguments at the hyper level.
%
\begin{thmm}
\label{thmm:sm-hm}
If $\rG$ has a structurally Markov graph law $\glaw$, and $\theta$
has a law from a hyper compatible set $\famlaw$, then
\[
(\rtheta_A, \rG_A) \indep(\rtheta_{B},
\rG_B) \bigmid \rtheta_{A \cap B}, \bigl\{\rG\in\decom{A} {B}\bigr\}\qquad [
\glaw, \famlaw].
\]
Furthermore, if each law $\mbox{\law}^{(\G)} \in\famlaw$ is strong hyper
Markov with respect to $\G$, then
\[
(\rtheta_A, \rG_A) \indep(\rtheta_{B|A},
\rG_B) \bigmid\bigl\{\rG\in \decom{A} {B}\bigr\}\qquad [\glaw, \famlaw].
\]
\end{thmm}
\begin{pf}
The proof for the first case is the same as in Theorem~\ref{thmm:sm-mark}. The proof for the strong case follows similar
steps, except starting with the strong hyper Markov property
\[
\rtheta_A \indep\rtheta_{B|A} \mid\rG, \bigl\{\rG\in\decom{A}
{B} \bigr\}. 
\]
\upqed
\end{pf}

Hyper compatible sets of strong hyper Markov laws have the additional
advantage that the posterior graph law will also be structurally
Markov: this follows from Theorem~\ref{thmm:sm-mark} and
\citet{dawid1993}, Proposition~5.6, which states that the marginal
distribution of the data under a strong hyper Markov law is Markov.
Furthermore, the posterior family of graph laws $\{ \mbox{\law}^{(\G)}(\cdot
\tmid X) \dvtx \G\in\sudg\}$ will maintain hyper compatibility.

\subsection{Clique vector}
\label{sec:clique-vector}

We show that the family of structural Markov laws forms an exponential
family of conjugate distributions for Bayesian updating under
compatible sampling.
%
\begin{dfn}
\label{dfn:charfn}
Define the \emph{completeness vector} of a graph to be the function
$c\dvtx \sudg\to\{0,1\}^{2^V}$ such that, for each $A \subseteq V$,
\[
c_A(\G) = \cases{ 1, & \quad$\mbox{if $\G_A$ is complete}$,
\vspace*{2pt}
\cr
0, &\quad $\mbox{otherwise}$. }
\]
Furthermore, define the \emph{clique vector} of a graph $t\dvtx \sudg\to
\Z^{2^V}$ to be the \mobius inverse of $c$ by \emph{superset}
inclusion
%
\begin{equation}
\label{eq:tcharfn} t_B(\G) = \sum_{A \supseteq B}
(-1)^{|A \setminus B|} c_A(\G).
\end{equation}
\end{dfn}

In the language of \citet{studeny2005}, $c$ and $t$ are both
\emph{imsets}.

The decomposition of $c$ and $t$ mirrors that of the graph.
%
\begin{lem}
\label{lem:ctdecom}
If $\G\in\decom{A}{B}$, then
%
\begin{eqnarray}
c(\G) &=& \bigl[c(\G_A)\bigr]^0 + \bigl[c(
\G_B)\bigr]^0 - \bigl[c(\G_{A \cap B})
\bigr]^0 \quad\mbox{and} \label{eq:cdecom}
\\
t(\G) &=& \bigl[t(\G_A)\bigr]^0 + \bigl[t(
\G_B)\bigr]^0 - \bigl[t(\G_{A \cap B})
\bigr]^0, \label{eq:tdecom}
\end{eqnarray}
where $[\cdot]^0$ denotes the expansion of a vector with zeroes to
the required coordinates.
\end{lem}
\begin{pf}
A subset $U \subseteq V$ induces a complete subgraph of $\G\in
\decom{A}{B}$ if and only if it induces a complete subgraph of
$\G_A$ or of $\G_B$ (or of both). \eqref{eq:cdecom} follows by the
inclusion-exclusion principle. \eqref{eq:tdecom} may then be
obtained by substitution into~\eqref{eq:tcharfn}.
\end{pf}
%
\begin{thmm}
\label{thmm:tform}
For any decomposable graph $\G\in\sudg$ and $A \subseteq V$,
\[
t_A(\G) = \cases{ 1, &\quad $\mbox{if } A \in\cl(\G)$, \vspace*{2pt}
\cr
-
\nu_\G(A), & \quad$ \mbox{if }A \in\csep(\G),\mbox{ and}$ \vspace *{2pt}
\cr
0, & \quad $\mbox{otherwise}$, }
\]
where $\cl(\G)$ are the cliques of $\G$, and $\csep(\G)$ are the
clique separators, and each separator $S$ has multiplicity
$\nu_\G(S)$.
\end{thmm}
\begin{pf}
For any $C \subseteq V$, let $\G^{(C)}$ be the graph on $V$ whose
edges are the set of all pairs $\{u,v\} \subseteq C$ (i.e.,
complete on $C$ and sparse elsewhere). Then it is straightforward
to see that
\[
t_A\bigl(\G^{(C)}_C\bigr) = \cases{ 1, &\quad $
\mbox{if } A=C$, \vspace*{2pt}
\cr
0, & \quad $\mbox{otherwise}$. }
\]
Now let $C_1,\ldots, C_k$ be a perfect ordering of the cliques of
$G$, and $S_2, \ldots, S_k$ be the corresponding separators. By
Lemma~\ref{lem:ctdecom}, it follows that
\[
t(\G) = \sum_{i=1}^k t\bigl(
\G^{(C_i)}_{C_i}\bigr) - \sum_{i=2}^k
t\bigl(\G ^{(S_i)}_{S_i}\bigr). 
\]
\upqed\end{pf}

Objects similar to the clique vector have arisen in several contexts.
Notably, it appears to be equivalent to the index $v$ of
Lauritzen, Speed and\break 
Vijayan [(\citeyear{lauritzen1984}), Definition~5],
which is characterised in a
combinatorial manner. It is also closely related to the
\emph{standard imset} of \citet{studeny2005}, which is equal to
\[
t\bigl(\G^{(V)}\bigr) - t(\G),
\]
where $\G^{(V)}$ is the complete graph.

The algorithm of \citet{wormald1985} for the enumeration of
decomposable graphs is based on a generating function for the vector
$\R^{|V|}$ that he termed the ``maximal clique vector,'' and is
equivalent to
\[
\operatorname{mcv}_k(\G) = \sum_{A \subseteq V \dvtx |A|=k}
t_A(\G), \qquad k=1,\ldots,|V|.
\]

\begin{prp}
\label{prp:tprop}
For any $\G\in\sudg$, the vector $t(\G)$ has the following
properties:
\begin{longlist}[(iii)]
\item[(i)]
\[
\sum_{A \subseteq V} t_A(\G) = 1 ,
\]
\item[(ii)] for each $v \in V$
\[
\sum_{A \ni v} t_A(\G) = 1,
\]
\item[(iii)]
\[
\sum_{A \subseteq V} |A| t_A(\G) = |V|\quad
\mbox{and}
\]
\item[(iv)]
\[
\sum_{A \subseteq V} \pmatrix{|A|
\cr
2} t_A(\G) =
\bigl|\edge(\G)\bigr|.
\]
\end{longlist}
\end{prp}
\begin{pf}
By the \mobius inversion theorem
[see, e.g., \citet{lauritzen1996}, Lemma~A.2], $c$ can also be expressed in terms of $t$,
\[
c_A(\G) = \sum_{B \supseteq A} t_B(
\G) ,\qquad  A \subseteq V.
\]
(i) and (ii) are $c_A(\G)$ at $A = \varnothing$ and $A = \{v\}$,
respectively, both of which induce complete subgraphs. (iii) is
obtained from (ii) by summation over $v \in V$, and (iv) is obtained
from (ii) by double counting each edge via summation over both
elements $\{u,v\} \in\edge(\G)$.
\end{pf}

\subsection{Clique exponential family}
\label{sec:exponential-family}

\begin{dfn}
\label{dfn:clique-expon}
The \emph{clique exponential family} is the exponential family of
graph laws over $\F\subseteq\sudg$, with $t$ as a natural
statistic (with respect to the uniform measure on $\sudg$). That
is, laws in the family have densities of the form
\[
\pi_\omega(\G) = \frac{1}{Z(\omega)} \exp\bigl\{ \omega\cdot t(\G)\bigr\}
, \qquad \G\in\F, \omega\in\R^{2^V},
\]
where $Z(\omega)$ is the normalisation constant, which will
generally be hard to compute.
\end{dfn}

Equivalently, the distribution can be parameterised in terms of $c$,
\[
\pi_\omega(\G) = \frac{1}{Z(\omega)} \exp \biggl\{ \biggl( \sum
_{B \subseteq A} (-1)^{|A \setminus B|} \omega_A
\biggr)_{A \subseteq V} \cdot c(\G) \biggr\},
\]
but $t$ is more useful due to the fact that it is sparse (by Theorem~\ref{thmm:tform}) and, as we shall see, is the natural statistic for
posterior updating.

Note that this distribution is over-parametrised.
By Proposition~\ref{prp:tprop}(i) and (ii), there are $|V|+1$ linear constraints in
the set of possible $t(G)$, adding multiples of $\alpha= (1)_{S
\subseteq V}$, or $\beta_v = (\ind_{v \in S})_{S \subseteq V}$ to
$\omega$ will leave the resulting $\pi$ unchanged. For the purpose of
identifiability, we could define a standardised vector $\omega^*$ as
\[
\omega^* = \omega+ \omega_\varnothing\alpha+ \sum
_{v \in V} (\omega _{\{v\}} - \omega_\varnothing)
\beta_v = \biggl( \omega_A+ \bigl(|A| - 1\bigr) \omega
_\varnothing- \sum_{v \in A} \omega_{\{v\}}
\biggr)_{A \subseteq V}
\]
such that $\pi_\omega= \pi_{\omega^*}$, and $\omega^*_{\{v\}} =
\omega^*_\varnothing= 0$ for all $v \in V$.

\begin{thmm}
\label{thmm:struct-markov-exponential}
Let $\glaw$ be a graph law whose support is $\sudg$. Then $\glaw$
is structurally Markov if and only if it is a member of the clique
exponential family.
\end{thmm}
\begin{pf}
See Appendix \ref{sec:proofs}.
\end{pf}

\begin{rem*}
It is possible to weaken the condition of full support; for example,
the same argument applies to any family $\F$ with the property that
if $\G\in\F$ and $C$ is a clique of $\G$, then $\G^{(C)} \in
\F$. 
\end{rem*}

A very similar family was proposed by \citet{bornn2011}; however, their
family allows the use of different parameters for cliques and
separators, which will generally not be structurally Markov.

\begin{exm}[{[{\citet{giudici1999}}; {\citet{brooks2003}, Section~8}]}]
\label{exm:uniform}
The simplest example of such a distribution is the uniform distribution over
$\sudg$, which by Proposition~\ref{prp:tprop}(i), corresponds to
$\omega_A$
being constant for all~$A$.
\end{exm}

\begin{exm}[{[{\citet{madigan1994,jones2005}}]}]
\label{exm:binomial}
Another common approach is to use a set of ${|V|\choose2}$
independent Bernoulli variables with probability $\psi$ to indicate
edge inclusion (i.e., an \erdosrenyi random graph), conditional on
$\rG$ being decomposable. The density of such a law is of the form
\[
\pi(\G) \propto\psi^{|\edge(\G)|} (1-\psi)^{{p\choose2} -
|\edge(\G)|} \propto \biggl(
\frac{\psi}{1-\psi} \biggr)^{|\edge(\G)|}.
\]
By Proposition~\ref{prp:tprop}(iv), it follows that this
distribution is a member of the exponential family with parameter
\[
\omega_A = \pmatrix{|A|
\cr
2} \log \biggl(\frac{\psi}{1-\psi
} \biggr).
\]
More generally, the family with parameter
\[
\omega_A = \sum_{e \in{A\choose2}} \log \biggl(
\frac{\psi_e}{1-\psi_e} \biggr)
\]
would correspond to the extension where each edge $e$ has its own
probability $\psi_e$.
\end{exm}
%
\begin{exm}
By adjusting parameters of the family, particular graphical features
can be emphasised. For example, a family of the form
\[
\omega_A = \pmatrix{|A|
\cr
2} \rho- \kappa\max\bigl(0,|A| - 2\bigr),
\]
with $\kappa> 0$, will penalise clique sizes greater than 2, placing a
higher probability on forest structures.
\end{exm}
%
\begin{exm}[{[\citet{armstrong2009}]}]
\label{exm:armstrong}
For comparison, it is useful to consider a nonstructurally Markov
graph law. Define the distribution over the number of edges to be
uniform, and the conditional distribution over the set of graphs
with a fixed number of edges to be uniform. This has density of the
form
\[
\pi(\G) = \frac{1}{{p\choose2} +1} \frac{1}{ | \{ \G' \in\sudg: |\edge(\G')|
= |\edge(\G)|  \} |}.
\]
Specifically for graphs on three vertices, we have that
\[
\pi\bigl(\mbox{\includegraphics{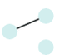}}\bigr) = \tfrac{1}{12},\qquad
\pi\bigl(\mbox{\includegraphics{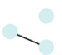}}\bigr) =
\tfrac{1}{12},\qquad \pi\bigl(\mbox{\includegraphics{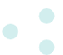}}\bigr) =
 \tfrac{1}{4}\quad
\mbox{and}\quad \pi\bigl(\mbox{\includegraphics{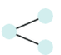}}\bigr)
= \tfrac{1}{12}.
\]
Therefore by Proposition~\ref{prp:sm-ratio} the law cannot be
structurally Markov.
\end{exm}

\subsection{Posterior updating}
\label{sec:posterior-updating}

We saw in Corollary~\ref{cor:sm-mark-post} that if the sampling
distributions are compatible, then posterior updating will preserve
the structural Markov property. In this section we show that this
updating may be performed locally, with the exponential clique family
forming a conjugate prior for a family of compatible models.

\begin{thmm}
\label{thmm:compatible-family-density}
Let $\famtheta$ be a family of compatible distributions for $X$,
where each $\theta^{(\G)}$ has density $\pi^{(\G)}$ with respect to
some product measure. Then
\[
\pi^{(\G)}(x) = \prod_{A \subseteq V} p_A
(x_A) ^{[t(\G)]_A},
\]
for all $x$ such that $\pi^{(\G)}(x) > 0$, where $p_A$ is the
marginal density of $X_A$ whenever $\G_A$ is complete, and
$p_\varnothing(x_\varnothing) = 1$.
\end{thmm}
\begin{pf}
For any decomposition $(A,B)$ of $\G$, then for any $x$ such that
$\pi^{(G)}(x) > 0$,
\begin{eqnarray*}
\pi^{(\G)}(x) &=& \pi_A^{(\G)}(x_A)
\pi_{B \mid A}^{(\G)}(x_{B
\setminus A} \mid x_A) =
\pi_A^{(\G)}(x_A) \pi_{B \mid A \cap
B}^{(\G)}(x_{B \setminus A}
\mid x_{A \cap B})
\\
&= &\pi_A^{(\G)}(x_A) \frac{\pi_B^{(\G)}(x_B)}{\pi_{A\cap
B}^{(\G)}(x_{A\cap B})}.
\end{eqnarray*}
The result follows by recursive decomposition over the clique tree.
\end{pf}

Therefore if the prior law for $\rG$ is a clique exponential with
parameter $\omega$, then under sampling from a compatible family the
resulting posterior law is of the same family,
\[
\pi(\G|X=x) \propto\exp \bigl\{ \bigl[\omega+ \bigl(\log p_A(x_A)
\bigr)_{A \subseteq V} \bigr] \cdot t(\G) \bigr\}.
\]

A key benefit of this conjugate formation is that we can describe the
posterior law with a parameter of dimension $2^{|V|}$ (strictly
speaking, we only need $2^{|V|} - |V| -1$, due to the
over-parametrisation). This is much smaller than for an arbitrary law
over the set of undirected decomposable graphs, which would require a
parameter of length approximately $2^{{|V|\choose2}}$.


\section{Ordered directed structural Markov property}
\label{sec:order-direct-struct}

We now investigate the first of two different methods by which the
structural Markov property might be extended to directed acyclic
graphical models (DAGs). In this section, we consider a law for a
random graph $\rdG$ over the set $\sodag$: the set of directed acyclic
graphs that respect a fixed well ordering $\prec$ on $V$.

The set $\sodag$ is straightforward to characterise, as $\prec$
determines the directionality of an edge between a pair of vertices.
Therefore, as in the undirected case, a random graph $\rdG$ on
$\sodag$ can also be interpreted as a random vector of length
${|V|\choose2}$.

In order to develop a structural Markov graph law over $\sodag$,
recall that the strong directed hyper Markov property can be expressed
as
%
\begin{equation}
\label{eq:sdhm-restate} \rtheta_{v \mid\pr(v)} \indep\rtheta_{\pr(v)} ,
\end{equation}
for all $v \in V$. This in turn implies mutual independence of the
collection $(\rtheta_{v \mid\pr(v)} )_{v \in V}$. Each element
$\rtheta_{v \mid\pr(v)}$ is constrained by $\G$ only through the
parent set $\pa_\G(v)$, as we require that $X_v \indep X_{\pr(v)}
\mid
X_{\pa_\G(v)}$. This motivates the following definitions.

\begin{dfn}
The \emph{ordered remainder graph of $\dG$ of $v \in V$ with respect
to $\prec$}, denoted by $\dG^\prec_{v \mid\pr(v)}$ is the graph
on $\{v\} \cup\pr(v)$, and edge set $\edge(\G_{\{v\} \cup\pr(v)})
\cup\{(w,u) \dvtx w,u \in\pr(v), w\prec u\}$, that is, the subgraph
$\G_{\{v\} \cup\pr(v)}$ with the addition of all possible edges
between elements of $\pr(v)$ respecting $\prec$.
\end{dfn}
The ordered remainder graph directly corresponds to the parent set of
the vertex, or equivalently, the set of vertices with a common head,
\[
\dG^\prec_{v \mid\pr(v)} \simeq\pa_\dG(v) \simeq\bigl\{
(u,w) \in \edge(\dG) \dvtx w = v\bigr\}.
\]
The advantage of the remainder graph is that it allows the
partitioning of the semi-graphoid into its constituent components.
%
\begin{prp}
\label{prp:order-partition}
Let $\dG$ be a directed acyclic graph compatible with the
ordering $\prec$. Then a distribution $P$ is Markov with respect to
$\dG$ if and only if for each $v \in V$, $P_{\{v\} \cup\pr(v)}$ is
Markov with respect to $\dG^\prec_{v \mid\pr(v)}$.

Similarly, a law $\mbox{\law}$ is weak/strong hyper Markov if and only if
for each $v \in V$, $\mbox{\law}_{\{v\} \cup\pr(v)}$ is weak/strong hyper
Markov with respect to $\dG^\prec_{v \mid\pr(v)}$.
\end{prp}
\begin{pf}
These follow from the ordered directed Markov property.
\end{pf}
The motivation of the term ``remainder'' is that $\dG^\prec_{v \mid
\pr(v)}$ encodes the remainder of the semi-graphoid of $\dG_{\{v\}
\cup\pr(v)}$ that is not determined by $\dG_{\pr(v)}$.

\begin{figure}

\includegraphics{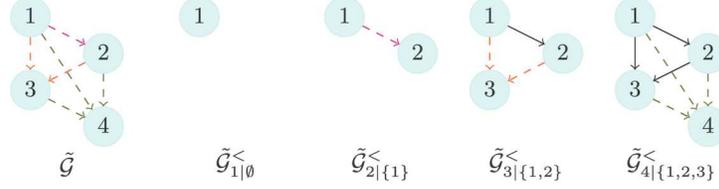}

\caption{A random directed acyclic graph $\rdG$ on $V =
\{1,2,3,4\}$, subject to the ordering $<$, and its corresponding
ordered remainder graphs. Under an ordered directed structural
Markov graph law, the ordered remainder graphs---or equivalently,
the collections of like-coloured edges---are independent.}
\label{fig:ordered-remainder}
\end{figure}

\begin{dfn}[(Ordered directed structural Markov property)]
The graph law $\glaw(\rdG)$ over $\sodag$ is \emph{ordered directed
structurally Markov} with respect to the ordering $\prec$ if for
each $v \in V$,
\[
\rdG^\prec_{v \mid\pr(v)} \indep\rdG_{\pr(v)}.
\]
\end{dfn}
As $\rdG_{\pr(v)} \simeq(\rdG^\prec_{u \mid\pr(v)})_{u \in\pr(v)}$,
this implies that the set of all ordered remainder graphs, or
equivalently, the set of all parent sets, are mutually independent; see
Figure~\ref{fig:ordered-remainder}.

Admittedly this construction is not very complicated, but it does
demonstrate how structure can be ``decomposed'' in directed graphs,
which will be used in the next section.


\section{Markov equivalence and the dagoid structural Markov property}
\label{sec:markov-equivalence}

The approach in Section~\ref{sec:order-direct-struct} cannot be
applied directly to distributions over the $\sdag$, the set of all
directed acyclic graphs on $V$. For instance, parent sets of
individual vertices cannot be independent: if $u$ is a parent of $v$,
then $v$ is precluded from being a parent of $u$.

A bigger problem is that there is no longer a one-to-one
correspondence between a graph and its semi-graphoid. That is, two or
more distinct DAGs may have identical conditional independence
properties, for example,
\mbox{\includegraphics{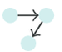}}, \mbox{\includegraphics[raise=-2pt]{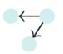}}, and
\mbox{\includegraphics{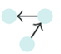}}.

\begin{dfn}
\label{dfn:dagoid}
Let $\dG$ and $\dG'$ be directed acyclic graphs such that
$\markov(\dG) = \markov(\dG')$. Then $\dG$ and $\dG'$ are termed
\emph{Markov equivalent}, and we write
\[
\dG\markoveq\dG'.
\]
A \emph{dagoid} is a Markov equivalence class of directed acyclic
graphs. We define the \emph{complete} and \emph{sparse} dagoids to
be the Markov equivalence classes of complete and sparse DAGs,
respectively. We use $\sdagoid$ to denote the set of dagoids on
$V$.
\end{dfn}

There are various methods of characterising Markov equivalence,
several of which are mentioned in the \hyperref[app]{Appendix}.

So when specifying a law for directed acyclic graphs, we are left with
the question of whether or not we should treat Markov equivalent
graphs as the same model. In other words, whether the model is
defined by the graph or the set of conditional independence statements
which it encodes. As noted earlier, we take the latter view.

A further advantage of working with equivalence classes is that a
smaller number of models needs be considered. Unfortunately this may
not be as beneficial as one may initially hope: \citet{castelo2004}
observed empirically that the ratio of the number DAGs to the number
of equivalence classes appears to converge to approximately $3.7$ as
the number of vertices increases.

\subsection{Ancestral sets and remainder dagoids}
\label{sec:prop-mark-equiv}

Although ancestral sets are used in the definition of the global
directed Markov property, ancestral sets themselves are not preserved
under Markov equivalence. However, as noted in
Section~\ref{sec:separoids}, subgraphs induced by ancestral sets
preserve the projection of the semi-graphoid. A somewhat trivial
consequence is the following.
%
\begin{prp}
\label{prp:dag-ancestral-equiv}
Let $\dG\markoveq\dG'$ and $A \subseteq V$ be ancestral in both
$\dG$ and $\dG'$. Then $\dG_A \markoveq\dG'_A$.
\end{prp}

This motivates the following definition.
%
\begin{dfn}
\label{dfn:dagoid-ancestral}
A set $A \subseteq V$ is \emph{ancestral} in a dagoid $\D$ if it is
ancestral for some graph $\dG\in\D$. For any such $A$, define the
\emph{subdagoid induced by $A$} to be the Markov equivalence class
of $\dG_A$, and denote it by $\D_A$.

For any $A \subseteq V$, let $\ances(A)$ denote the set of dagoids
on $V$ in which $A$ is an ancestral set.
\end{dfn}
Note that the dagoid ancestral property is not as strong as the
collapsibility property in undirected graphs, in that there can exist
nonancestral sets that also preserve the semi-graphoid of the induced
subgraph.

However, ancestral sets are still quite powerful, in that they can be
used to decompose the semi-graphoid.
%
\begin{dfn}
\label{dfn:dag-insertion}
Let $\dG$ be a directed acyclic graph on $V$, of which $A$ is an
ancestral set, and let $\dgraph{H}$ be a directed acyclic graph on
$A$. Then the \emph{insertion of $\dgraph{H}$ into $\dG$}, written
\[
\dgraph{H} \dprod\dG,
\]
is the directed acyclic graph on $V$ with edge set
\[
\edge(\dgraph{H}) \cup \bigl[ \edge(\dG) \setminus A^2 \bigr].
\]
\end{dfn}
In other words, the edges between elements of $A$ are determined by
$\dgraph{H}$, and all other edges are determined by $\dG$. This
operation preserves Markov equivalence.
%
\begin{lem}
\label{lem:dag-insertion-equiv}
Let $\dG$ and $\dG'$ be Markov equivalent graphs in which $A$ is an
ancestral set, and $\dgraph{H}$ and $\dgraph{H}'$ be Markov
equivalent graphs on $A$. Then
\[
\dgraph{H} \dprod\dG\markoveq\dgraph{H}' \dprod\dG'.
\]
\end{lem}
\begin{pf}
We use the notation and results of Appendix~\ref{sec:mark-equiv}.
Both graphs must have the same skeleton. Let $(a,b,c)$ be an
immorality in $\dgraph{H} \dprod\dG$. Then if $b \in A$, then
$(a,b,c)$ must be an immorality of $\dgraph{H}$, and hence also an
immorality of $\dgraph{H}'$, and so also of $\dgraph{H}' \dprod
\dG'$.

Otherwise if $b \notin A$, and at least one of $a$ or $c$ is not in
$A$, then $(a,b,c)$ must be an immorality of $\dG$, and hence an
immorality of $\dG'$ and $\dgraph{H}' \dprod\dG'$.

Finally, if $b \notin A$ and $a,c \in A$, then $\{a,c\}$ must not be
an edge in the skeleton $\dgraph{H}$, nor an edge in the skeleton of
$\dgraph{H}'$. Hence it must also be an immorality of $\dgraph{H}'
\dprod\dG'$.
\end{pf}

Consequently for a dagoid $\D$ with ancestral set $A$, we can define
the \emph{ancestral insertion} of a dagoid $\dagoid{K}$ on $A$ into
$\D$ as
\[
\dagoid{K} \dprod\D= [\dgraph{H} \dprod\dG],
\]
where $\dG\in\D$ is a directed acyclic graph with an ancestral set
$A$, $\dgraph{H} \in\dagoid{K}$, and $[\cdot]$ denotes the Markov
equivalence class.

\begin{figure}

\includegraphics{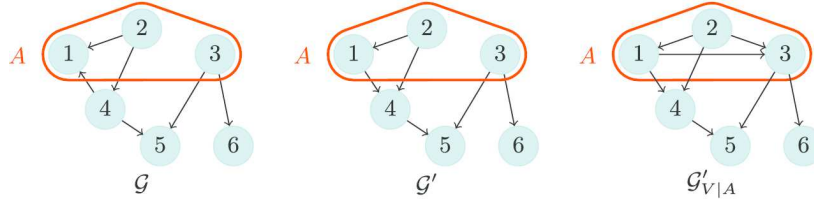}

\caption{$A = \{1,2,3\}$ is an ancestral set of the dagoid
containing $\dG$, as it is ancestral in the graph $\dG'$ obtained
by reversing the covered edge $(4,1)$. $\dG'_{V \mid A}$ is
obtained by replacing the edges between elements of $A$ with those
of a complete graph on $A$.}
\label{fig:dagoids}
\end{figure}

We use this approach to extend the notion of a remainder graph from
the previous section without the use of a fixed well-ordering.
%
\begin{dfn}
\label{dfn:dag-remainder}
Let $A$ be an ancestral set of a directed acyclic graph $\dG$.
A~directed acyclic graph $\dG_{V|A}$ is a \emph{remainder graph of
$\dG$ given $A$} if
\[
\dG_{V|A} = \dagoid{C}^{(A)} \dprod\dG,
\]
where $\dagoid{C}^{(A)}$ is a complete dagoid on $A$.

By Lemma~\ref{lem:dag-insertion-equiv}, the remainder graph must be
unique up to Markov equivalence. Hence for a dagoid $\D\in
\ances(A)$, we can uniquely define the \emph{remainder dagoid of
$\D$ given $A$}, denoted by $\D_{V|A}$; see
Figure~\ref{fig:dagoids}.
\end{dfn}

Analogous with the ordered case, the induced and remainder dagoids
$\D_A$ and $\D_{V \mid A}$ characterise the complete dagoid (via the
ancestral insertion). Moreover, they can be chosen independently.
%
\begin{thmm}
\label{thmm:dagoid-variation-indep}
For any $A \subseteq V$, we have
\[
\D_A \vindep\D_{V|A} \mid\bigl\{ \D\in\ances(A)\bigr\}.
\]
\end{thmm}
\begin{pf}
For any $\D,\D' \in\ances(A)$, we can construct $\D^* = \D_A
\dprod
\D'_{V|A}$. This will have the required properties that $\D^*_A =
\D_A$ and $\D^*_{V|A} = \D'_{V|A}$.
\end{pf}


\subsection{Dagoid structural Markov property}
\label{sec:dag-smp}

This motivates the following construction for the structural Markov
property.
%
\begin{dfn}[(Dagoid structural Markov property)]
\label{dfn:dag-smp}
We say a graph law $\glaw(\rD)$ is \emph{structurally Markov} if for
any $A \subseteq V$, we have
\[
\rD_{V|A} \indep\rD_A \mid\bigl\{\rD\in\ances(A)\bigr\} \qquad [\glaw].
\]
\end{dfn}

As in the undirected case, we can characterise this property via the
odds ratio of the density.
%
\begin{prp}
\label{prp:dag-sm-ratio}
A graph law is structurally Markov if and only if for any $\D,\D'
\in\ances(A)$, we have
%
\begin{equation}
\label{eq:dag-smp-crossprod} \pi(\D) \pi\bigl(\D'\bigr) = \pi\bigl(
\D_A \dprod\D'_{V|A}\bigr) \pi\bigl(
\D'_A \dprod \D_{V|A}\bigr).
\end{equation}
\end{prp}
\begin{pf}
As in Proposition~\ref{prp:sm-ratio}, we may write the density
\[
\pi\bigl(\D\tmid\ances(A)\bigr) = \pi\bigl(\D_A \tmid\ances(A)
\bigr) \pi\bigl(\D_{V|A} \tmid\ances(A)\bigr).
\]
\upqed\end{pf}

\subsection{d-Clique vector}
\label{sec:d-clique-vector}

The equivalence class formulation of a dagoid is difficult to work
with, both algebraically and computationally. Instead we propose a
characteristic vector similar to the clique vector of Section~\ref{sec:clique-vector}.
%
\begin{dfn}
\label{dfn:d-clique}
The \emph{d-clique vector} of a directed acyclic graph $\dG$ is
%
\begin{equation}
\label{eq:d-clique} t(\dG) = \sum_{v \in V} \bigl[\delta
\bigl(\{v\} \cup\pa_{\dG} (v)\bigr) - \delta\bigl(\pa_{\dG}
(v)\bigr) \bigr] + \delta(\varnothing) \in\Z^{2^V},
\end{equation}
where $\delta(A) = (\ind_{S = A})_{S \subseteq V}$.
\end{dfn}
Again, we note the relationship to the imsets of \citet{studeny2005},
specifically the structural imset $u_{\dG} = \delta(V) - t(\dG)$ in
Section~\ref{sec:mark-equiv-imsets}. For our purposes, the d-clique
vector is a more convenient object with which to work. This exhibits
analogous properties to those of the clique vector of
Section~\ref{sec:clique-vector}.

\begin{prp}
\label{prp:dagoid-tprop}
The properties of Proposition~\ref{prp:tprop} apply to all directed
graphs $\dG\in\sdag$.
\end{prp}
\begin{pf}
(i) follows directly from the definition. (ii) is obtained by
noting that each term of \eqref{eq:d-clique} contributes 1 if the
summand is $v$, and 0 otherwise. For (iii), each term of
\eqref{eq:d-clique} contributes 1, and (iv) is due to each term of
\eqref{eq:d-clique} counting the number of edges whose head is $v$.
\end{pf}

In a similar manner to the undirected case, we can define the
\emph{d-completeness vector} to be the \mobius transform of the
d-clique vector,
%
\begin{equation}
\label{eq:d-complete} c_A(\dG) = \sum_{B \supseteq A}
t_B(\dG),
\end{equation}
and say that a set $A \subseteq B$ is \emph{d-complete} if $c_A(\dG) =
1$. This corresponds to the definition of the characteristic imset of
\citet{hemmecke2012}.

\begin{lem}[({\citet{hemmecke2012}, Theorem~1})]
\label{lem:d-complete-parents}
Let $\prec$ be a well-ordering of a directed acyclic graph $\dG$.
For any nonempty set $A \subseteq V$, with maximal element $a$
under $\prec$,
\[
c_A (\dG) = \cases{ 1, &\quad $\mbox{if } A \setminus\{a\} \subseteq
\pa_{\dG} (a) $, \vspace *{2pt}
\cr
0, & \quad $\mbox{otherwise.}$ }
\]
\end{lem}

This provides the link to the completeness and clique vectors of
undirected graphs from Section~\ref{sec:clique-vector}.
%
\begin{cor}
\label{cor:d-clique-perfect}
If $\dG$ is a perfect directed acyclic graph, and $\dG^{\mathrm{s}}$
is its skeleton, then $c_{\dG} = c_{\dG^{\mathrm{s}}}$, and hence
$t(\dG) = t(\dG^{\mathrm{s}})$.
\end{cor}

Most important, the d-clique vector is a unique representation of the
dagoid.
%
\begin{thmm}
\label{thmm:d-clique-markov-equiv}
Let $\dG,\dG'$ be directed acyclic graphs on $V$. Then $\dG
\markoveq\dG'$ if and only if $t(\dG) = t(\dG')$.
\end{thmm}
\begin{pf}
To show that the d-clique vector is preserved under Markov
equivalence, by Theorem~\ref{thmm:dag-equiv-reversal} it is
sufficient to show that it is preserved under a covered edge
reversal. If $(a,b)$ is a covered edge of $\dG$, then the
contribution of these vertices to the sum \eqref{eq:d-clique} is
\begin{eqnarray*}
t(\dG) &=& \bigl[\delta\bigl(\{a\} \cup\pa_\dG(a)\bigr) - \delta
\bigl(\pa_\dG (a)\bigr) \bigr] + \bigl[\delta\bigl(\{b\} \cup
\pa_\dG(b)\bigr) - \delta\bigl(\pa_\dG(b)\bigr) \bigr]
\\
& &{}+ \sum_{v \ne a,b} \bigl[ \delta\bigl(\{b\} \cup
\pa_\dG(b)\bigr) - \delta\bigl(\pa_\dG(b)\bigr) \bigr] +
\delta(\varnothing).
\end{eqnarray*}
By definition, $\pa_{\dG}(a) \cup\{a\} = \pa_{\dG}(b)$, and so the
corresponding terms will cancel. If $\dG^*$ is obtained from $\dG$
by reversing $(a,b)$, note that
\[
\pa_{\dG}(a) = \pa_{\dG^*}(b) \quad\qquand\quad  \pa_{\dG}(b)
\cup\{b\} = \pa_{\dG^*}(a) \cup\{a\},
\]
and the remaining terms will be unchanged. Hence $t(\dG) =
t(\dG^*)$.

To show that the d-completeness vector (and hence, also the d-clique
vector) is unique to the equivalence class, by Theorem~\ref{thmm:dag-equiv-immoral} we can show that it determines the
skeleton and immoralities. By Lemma~\ref{lem:d-complete-parents},
there is an edge between $u$ and $v$ in $\dG$ if and only if
$c_{\{u,v\}}(\dG) = 1$. Likewise, $(u,v,w)$ is an immorality if and
only if $c_{\{u,v,w\}}(\dG) = 1$ and $c_{\{u,w\}}(\dG) = 0$.
\end{pf}

\begin{figure}

\includegraphics{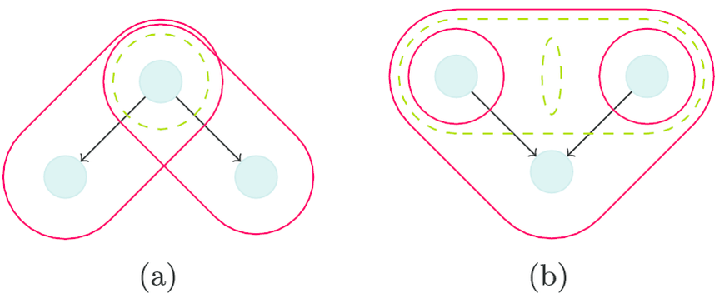}

\caption{The d-cliques (\protect\includegraphics{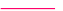}) and d-separators
(\protect\includegraphics{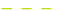}) of different directed acyclic graphs. Note
that in the perfect DAG \textup{(a)}, the d-cliques and d-separators are
the cliques and separators of the skeleton. However, as in \textup{(b)},
d-separators may contain d-cliques. %
}
\label{fig:d-clique}
\end{figure}

This cancellation of terms involving covered edges is very useful: as
a consequence, the d-clique vector will generally be quite sparse. In
line with the clique vector, we term a set $A \subseteq V$ such that
$t_A(\D)=1$ a \emph{d-clique}, and set $A$ such that $t_A(\D)<0$ a
\emph{d-separator}: See examples in Figure~\ref{fig:d-clique}.

%
%
%
%
%
\begin{thmm}
\label{thmm:d-clique-ancestral}
Let $A$ be an ancestral set of a dagoid $\D$. Then
\[
t(\D) = \bigl[t(\D_A)\bigr]^0 + t(\D_{V|A}) -
\delta(A),
\]
where $[\cdot]^0$ denotes the expansion of the vector with zeroes to
the required coordinates.
\end{thmm}
\begin{pf}
Let $\dG\in\D$ in which $A$ is ancestral, and $\prec$ be a
well-ordering of $\dG$ in which elements of $A$ precede those of $V
\setminus A$. Then
\[
\pa_{\dG} (v) = \cases{ \pa_{\dG_A} (v), & \quad $v \in A$,
\vspace*{2pt}
\cr
\pa_{\dG_{V|A}} (v), & \quad $v \notin A$. }
\]
The result follows after noting that
\[
\sum_{v \in A} \bigl[ \delta\bigl(\pa_{\dG_{V|A}}(v)
\cup\{v\}\bigr) - \delta\bigl(\pa_{\dG
_{V|A}}(v)\bigr) \bigr] = \delta(A).
\]
\upqed\end{pf}

We now arrive at the key result of this section: the dagoid structural
Markov property characterises an exponential family of graph laws.
%
\begin{thmm}
\label{thmm:dsmp-exp}
Let $\glaw$ be a graph law whose support is $\sdagoid$. Then
$\glaw$ is structurally Markov if and only if it is a member of the
exponential family with the d-clique vector as natural sufficient
statistic, that is, if $\glaw$ has density of the form
%
\begin{equation}
\label{eq:dsmp-exp} \pi_\omega(\D) \propto\exp\bigl\{ \omega\cdot t(\D)
\bigr\}.
\end{equation}
\end{thmm}
\begin{pf}
See Appendix~\ref{sec:proofs}.
\end{pf}

\begin{exm}
\label{exm:dagoid-uniform}
As in the undirected case, the simplest example of a structurally
Markov graph law is the uniform law over $\sdagoid$, on taking
$\omega_A = 0$.
\end{exm}
%
\begin{exm}
\label{exm:dagoid-erdos}
For any directed graph $\dG\in\D$, let $e(\D)$ denote
$|\edge(\dG)|$. By Proposition~\ref{prp:dagoid-tprop}, $e(\D) =
\sum_A {|A|\choose2} t_A(\D)$. So, for any $\rho> 0$, the graph
law specified by
\[
\pi(\D) \propto\rho^{e(\D)}
\]
is structurally Markov, on taking $\omega_A = {|A|\choose2} \log
\rho$.
\end{exm}

However, we note that some simple laws are \emph{not} structurally
Markov.
%
\begin{exm}
\label{exm:dag-uniform}
Consider the law in which $\pi(\D)$ is proportional to $|\D|$, in
other words, the uniform law on $\sdag$ projected onto $\sedag$.
Then using $[\cdot]$ to denote Markov equivalence class, we note the
size of the following dagoids:
\begin{eqnarray*}
\bigl[\mbox{\includegraphics{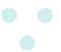}}\bigr] &=&
\bigl\{\mbox{\includegraphics{1319i12.eps}} \bigr\},
\\
{}\bigl[\mbox{\includegraphics{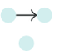}}\bigr]
&=& \bigl\{\mbox{\includegraphics{1319i13.eps}},\mbox{\includegraphics[raise=-3pt]{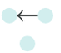}}\bigr\},
\\
{}\bigl[\mbox{\includegraphics{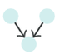}}\bigr] &=&
\bigl\{\mbox{\includegraphics{1319i15.eps}}\bigr\},
\\
{}\bigl[\mbox{\includegraphics{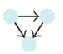}}\bigr] &=&
 \bigl\{\mbox{\includegraphics{1319i16.eps}},
\mbox{\includegraphics{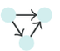}}, \mbox{\includegraphics[raise=-3pt]{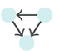}},
 \mbox{\includegraphics{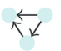}},
\mbox{\includegraphics{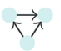}},
 \mbox{\includegraphics{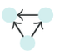}}\bigr\}.
\end{eqnarray*}
As a consequence, this law does not satisfy the property
$\pi\bigl(\bigl[\mbox{\includegraphics{1319i16.eps}}\bigr]\bigr) 
\pi\bigl(\bigl[\mbox{\includegraphics[raise=-3pt]{1319i12}}\bigr]\bigr) =
\pi\bigl(\bigl[\mbox{\includegraphics{1319i13.eps}}\bigr]\bigr) 
\pi\bigl(\bigl[\mbox{\includegraphics[raise=-3pt]{1319i15}}\bigr]\bigr)$ required by
Proposition~\ref{prp:dag-sm-ratio}.
\end{exm}

We note that similar exponential families were proposed by
\citet{mukherjee2008}. However, they treat Markov equivalent graphs
as distinct, and allow them to have different probabilities.

\subsection{Compatible distributions and laws}
\label{sec:dagoid-compatible}

As with the undirected case, a graph law is only part of the story.
For each dagoid $\D$, we also require a method to specify a Markov
sampling distribution and a law over such sampling distributions.

\begin{dfn}
\label{dfn:dag-compat}
Distributions $\theta$ and $\theta'$, Markov with respect to
directed acyclic graphs $\dG$ and $\dG'$ respectively, are termed
\emph{graph compatible} if, for every vertex $v$ such that
$\pa_{\dG}(v) = \pa_{\dG'}(v)$, there exist versions of the
conditional probability distributions for $X_v \mid X_{\pa(v)}$ such
that
\[
\theta(X_v \tmid X_{\pa(v)}) = \theta'(X_v
\tmid X_{\pa(v)}).
\]
Distributions $\theta$ and $\theta'$, Markov with respect to dagoids
$\D$ and $\D'$, respectively, are termed \emph{\textup{(}dagoid\textup{)} compatible}
if they are graph compatible for every pair of graphs $\dG\in\D,
\dG' \in\D'$.

Likewise, laws $\mbox{\law}(\rtheta)$ and $\mbox{\law}'(\rtheta)$, hyper Markov
with respect to $\dG$ and $\dG'$, respectively, are termed
\emph{graph hyper compatible} if for every vertex $v$ such that
$\pa_{\dG}(v) = \pa_{\dG'}(v)$, there exist versions of the
conditional laws for $\rtheta_{v|\pa(v)} \mid\rtheta_{\pa(v)}$ such
that
\[
\mbox{\law}(\rtheta_{v|\pa(v)} \tmid\rtheta_{\pa(v)}) = \mbox{\law}'(
\rtheta_{v|\pa(v)} \tmid\rtheta_{\pa(v)}).
\]
By \citet{dawid2001}, Section~8.2, the weak hyper Markov property
may be characterised in terms of $\markov(\dG)$, and so the weak
hyper Markov property can be defined with respect to a dagoid. Laws
$\mbox{\law}(\rtheta)$ and $\mbox{\law}'(\rtheta)$, that are hyper Markov with
respect to $\D$ and $\D'$, respectively, are \emph{\textup{(}dagoid\textup{)} hyper
compatible} if they are graph compatible for every pair of graphs
$\dG\in\D, \dG' \in\D'$.
\end{dfn}

As in the undirected case, we can define a family of compatible
distributions $\famtheta=\{\theta^{(\G)} \dvtx \G\in\sudg\}$ and a
family of hyper compatible laws $\famlaw= \{\mbox{\law}^{(\G)} \dvtx \G\in
\sudg\}$ if they are pairwise compatible or hyper compatible with
respect to the relevant graphs.

\begin{prp}
\label{prp:dagoid-compat-ci}
Suppose $\glaw(\rD)$ is a graph law over $\sdagoid$ and $\famtheta$
is a family of compatible distributions. Then
%
\begin{equation}
\label{eq:dagoid-compat-ci-1} X_A \indep\rD_{V|A} \mid
\rD_A, \bigl\{\rD\in\ances(A)\bigr\} \qquad[\famtheta,\glaw]
\end{equation}
and
%
\begin{equation}
\label{eq:dagoid-compat-ci-2} X_{V \setminus A} \indep\rD_A \mid X_A,
\rD_{V|A}, \bigl\{\rD\in\ances (A)\bigr\} \qquad[\famtheta,\glaw].
\end{equation}
Likewise, if $\glaw(\rD)$ is a graph law over $\sdagoid$ and
$\famlaw$ is a hyper compatible family of laws, then
\[
\rtheta_A \indep\rD_{V|A} \mid\rD_A, \bigl\{
\rD\in\ances(A)\bigr\}\qquad [\famlaw,\glaw]
\]
and
\[
\rtheta_{V\setminus A|A} \indep\rD_A \mid\rtheta_A,
\rD_{V|A}, \bigl\{ \rD\in\ances(A)\bigr\}\qquad [\famlaw,\glaw].
\]
\end{prp}
\begin{pf}
This is much the same as Proposition~\ref{prp:compat-ci}: for
\eqref{eq:dagoid-compat-ci-1}, the distribution of $X_A$ is
determined by the parent sets of the vertices in $A$ in some $\G\in
\D$ in which $A$ is ancestral. Likewise, in
\eqref{eq:dagoid-compat-ci-2}, the conditional distribution for
$X_{V \setminus A} \mid X_A$ is determined by the parent sets of
vertices in $V \setminus A$. The same argument applies at the hyper
level.
\end{pf}

Note that in the definition of compatibility and hyper compatibility
we specifically refer to \emph{versions} of conditional probabilities
and laws, as in some cases the conditional distributions/laws will not
be uniquely defined, due to conditioning on null sets.

As the weak hyper Markov property is defined on the separoid, the weak
directed hyper Markov property is well-defined for any dagoid.
However, the strong form requires further conditions.
%
\begin{dfn}
\label{dfn:strong-dagoid-hyper-Markov}
A law $\mbox{\law}(\rtheta)$ over $\markovdist(\D)$ is \emph{strong hyper
Markov} with respect to $\D$ if it is strong directed hyper Markov
with respect to every $\dG\in\D$.
\end{dfn}

If $\dG\in\D$ is perfect, then the strong dagoid hyper Markov
property is equivalent to the undirected strong hyper Markov property
on the skeleton of $\dG$; see \citeauthor{dawid1993}
[(\citeyear{dawid1993}), Proposition~3.15].
The notion of hyper compatibility is equivalent to the ``parameter
modularity'' property of \citet{heckerman1995}. Likewise, the strong
hyper Markov property is equivalent to their ``parameter
independence.''

\begin{exm}
\label{exm:dagoid-hiw}
For each vertex $v$ of a directed acyclic graph $\dG$, we define the
law for the conditional parameter $\mbox{\law}(\rtheta_{v|\pa_{\dG}(v)})$
to be the same as that of the inverse Wishart $\dinvwish(\nu;
\Phi)$. That is, using the notation of \citet{dawid1981}, we have
\[
\theta_{v | \pa_{\dG}(v)} = \dnorm(\Gamma_{v | \pa_\dG(v)}, \Sigma_{v | \pa_\dG(v)}),
\]
where
\begin{eqnarray*}
\mbox{\law}(\rSigma_{v|\pa_{\dG}(v)}) &=& \dinvwish\bigl(\nu+
\bigl|\pa_{\dG}(v)\bigr|;
\Phi_{v|\pa_{\dG}(v)}\bigr), 
\\
\mbox{\law}(\rGamma_{v|\pa_{\dG}(v)} \tmid\rSigma_{v|\pa_{\dG}(v)}) &=&
\Phi_{\{v\},\pa_{\dG}(v)}\Phi_{\pa_{\dG}(v)}^{-1} + \dnorm_{\{v\}
\times\pa_{\dG}(v)}
\bigl(\rSigma_{v|\pa_{\dG}(v)}, \Phi_{\pa_{\dG}(v)}^{-1}\bigr).
\end{eqnarray*}
By the properties of the inverse Wishart law, it follows that the
law is preserved under covered edge reversals. Therefore by Theorem~\ref{thmm:dag-equiv-reversal}, it is well defined for a dagoid, and
so may be termed the \emph{dagoid hyper inverse Wishart law}. Note
that this property is not satisfied by the more general inverse
type-II Wishart family of \citet{letac2007}.
\end{exm}

\begin{thmm}
\label{thmm:dagoid-struct-markov-strong-compat}
If $\famlaw$ is a family of strong hyper Markov hyper compatible
laws, then the family of marginal data distributions is compatible.
\end{thmm}
\begin{pf}
The hyper compatibility and the strong hyper Markov property imply
that, for any two dagoids $\D,\D'$ and any $\dG\in\D, \dG' \in
\D'$, if $\pa_{\dG}(v) = \pa_{\dG'}(v)$ for some $v \in V$, then
\[
\mbox{\law}^{(\D)}(\rtheta_{v|\pa}) = \mbox{\law}^{(\D')}(
\rtheta_{v|\pa}).
\]
Therefore, the family of marginal data distributions $\bar\famtheta
= \{\bar\theta^{(\D)}\dvtx \D\in\sdagoid\}$ will have
\[
\bar\theta^{(\D)}(X_v \tmid X_{\pa_{\dG}}) =
\E_{\mbox{\scriptsize{\law}}}^{(\D)} [ \rtheta_{v|\pa_{\dG}}] = \bar\theta^{(\D')}(X_v
\tmid X_{\pa_{\dG}}) = \E_{\mbox{\scriptsize{\law}}}^{(\D')} [ \rtheta_{v|\pa_{\dG}}].
\]
\upqed\end{pf}

This is particularly useful because, as in the undirected case, the
structural Markov property will be preserved in the posterior under
compatible sampling.
%
\begin{thmm}
\label{thmm:dagoid-struct-markov-post}
Suppose $\glaw(\rD)$ is a structurally Markov graph law over
$\sdagoid$ and $\famtheta$ is a family of compatible distributions.
Then the posterior graph law for $\rD$ is structurally Markov.
\end{thmm}
\begin{pf}
By the structural Markov property and \eqref{eq:dagoid-compat-ci-1},
we have
\[
(X_A,\rD_A) \indep\rD_{V|A} \bigmid\bigl\{
\rD\in\ances(A)\bigr\} ,
\]
and hence
\[
\rD_A \indep\rD_{V|A} \bigmid X_A, \bigl\{
\rD\in\ances(A)\bigr\}.
\]
Combining this with \eqref{eq:dagoid-compat-ci-2}, we get
\[
\rD_A \indep(\rD_{V|A},X_{V \setminus A}) \bigmid
X_A, \bigl\{\rD\in \ances(A)\bigr\} ,
\]
and hence
\[
\rD_A \indep\rD_{V|A} \bigmid X, \bigl\{\rD\in\ances(A)
\bigr\}. 
\]
\upqed\end{pf}

\subsection{Posterior updating}
\label{sec:posterior-updating-1}

If it is possible to avoid the problem of conditioning on null sets,
then as in the undirected case, a compatible family can be
characterised by a distribution on the complete dagoid.
%
\begin{thmm}
\label{thmm:dag-compat-complete}
If the distribution on the complete dagoid has positive density $p$
with respect to some product measure, then the compatible
distribution for any dagoid $\D$ has density
%
\begin{equation}
\label{eq:dag-compat-complete} \pi^{(\D)}(x) = \prod_{A \subseteq V}
p(x_A)^{[t(\D)]_A}.
\end{equation}
\end{thmm}
\begin{pf}
Let $\dG$ be an arbitrary graph in $\D$. Then by compatibility,
\[
p^{(\D)}(x) = \prod_{v \in V} p(x_v
| x_{\pa(v)}) = \frac{
\prod_{i = 1}^p p(x_{\{v_i\} \cup\pa(v_i)})
}{
\prod_{i = 2}^p p(x_{\pa(v_i)})
} = \prod_{A \subseteq V}
\bigl[p(x_A)\bigr]^{t(\D)_A}. 
\]
\upqed\end{pf}

As a consequence, if the graph law has a d-clique exponential family
of the form~\eqref{eq:dsmp-exp}, and the sampling distributions are
compatible with density of the form~\eqref{eq:dag-compat-complete},
then the posterior graph law will have density
\[
\pi(\D\tmid X) \propto\exp \bigl\{ \bigl[ \omega+ \bigl(\log
p_A(X_A)\bigr)_{A \subseteq V} \bigr] \cdot t(\D) \bigr
\}.
\]
That is, the d-clique exponential family is a conjugate prior under
sampling from a compatible family.


\section{Discussion}
\label{sec:discussion}

We have demonstrated how conditional independence can be used to
characterise families of distributions over undirected graphs, ordered
directed graphs, and equivalence classes of directed graphs.

One point to emphasise is that all three structural Markov properties
are distinct, in that no one property can be derived as a special case
of another; for example, the undirected structural Markov property does not
arise from the dagoid structural Markov property restricted to
equivalence classes of perfect DAGs.

\subsection{Open questions}

One significant open question is how the full support requirements of
Theorems~\ref{thmm:struct-markov-exponential} and \ref{thmm:dsmp-exp}
might be weakened. Obviously these theorems would not hold for all
subsets of graphs/dagoids, though we conjecture that they will hold
for any structurally meta Markov subsets.
A related problem is characterising structurally meta Markov subsets
of graphs.

\subsection{Computation}

One problem which we have not broached is the numerical calculation of
such graph laws. Except for the ordered directed case, where the
computations can be done in parallel, for even small numbers of
vertices it can quickly become infeasible to enumerate all graphs, and
hence some sort of numerical approximation will usually be required.
Markov chain Monte Carlo (MCMC) methods are commonly utilised for this
purpose.

For undirected decomposable graphs, \citet{giudici1999} proposed a
method in
which each iteration proposes adding or removing a single edge. They consider
the problem of sampling from the posterior of a uniform prior with a
compatible family sampling distributions, though this procedure can be applied
to any structural Markov graph law. This requires computing the
Metropolis--Hastings acceptance ratio,
\[
\min \biggl( \frac{\pi(\G')}{\pi(\G)},1 \biggr) =
\cases{ \min \bigl( \exp\bigl\{
\omega\cdot \bigl[t\bigl(\G'\bigr) - t(\G)\bigr]\bigr\},1 \bigr), & \quad$
\G,\G' \in\sudg$, \vspace *{2pt}
\cr
0, & \quad $\mbox{otherwise.}$ }
\]
The results of \citeauthor{frydenberg1989} [(\citeyear{frydenberg1989}),
Lemma~3] and Giudici and\break Green [(\citeyear{giudici1999}),
Theorem~2] characterise such so-called neighbouring graphs, and
also imply that for any two such graphs $\G, \G'$, the vector $t(\G')
- t(\G)$ has only 4 nonzero elements. Consequently, for any
structurally Markov graph law over $\sudg$, the parameter $\omega$
need only be evaluated on 4 such places: this is particularly
beneficial for posterior graph laws where each element of $\omega$
requires the evaluation of the marginal density of the model.

Unfortunately, such algorithms often exhibit poor mixing properties
[\citet{kijima2008}], resulting in unreliable estimates.
\citet{green2013} develop an extension for making proposals which add
or remove multiple edges, resulting in faster mixing: this algorithm
is also able to take advantage of local computations in computing the
acceptance ratio.

For dagoids, the problem is considerably more difficult.
\citet{chickering2003}, \citet{auvray2002} and \citet{studeny2005a}
have developed methods for characterising the neighbouring dagoids
(i.e., dagoids obtained by adding or removing an edge to a graph in
the current dagoid). \citet{he2013} recently developed an MCMC scheme
based on
this approach: as in the undirected case, the acceptance ratio will
also depend on a sparse vector, and so can be computed efficiently.

More generally, the problem of finding the most probable graph under a
structural Markov law, which for posterior laws is known as maximum
\textit{a posteriori} (MAP) estimation, is an example of a
\emph{strong decomposable search criterion}
[\citet{studeny2005}, Section~8.2.3]. As suggested by \citet{hemmecke2012}, linear and
integer programming techniques based on the (d-)clique or
(d-)completeness vectors may provide elegant solutions to this problem.

\subsection{Extensions}

A further open question is how structural Markov properties might be
defined for other classes of graphical models, such as
nondecomposable undirected graphs, ancestral graphs, and marginal
independence (bidirected) graphs. The
identification of such properties would rely on establishing
constructions for partitioning the structure, analogous to
decompositions and ancestral graphs.


\begin{appendix}\label{app}
\section{Characterising Markov equivalence of directed acyclic graphs}
\label{sec:mark-equiv}

Numerous techniques have been developed for determining whether two
graphs are Markov equivalent.

\subsection{Skeleton and immoralities}
\label{sec:skel-immor}

The \emph{skeleton} of a DAG is the undirected graph obtained by
substituting the directed edges for undirected ones. A triplet
$(a,b,c)$ of vertices is an \emph{immorality} of a DAG $\dG$ if the
induced graph $\dG_{\{a,b,c\}}$ is of the form $a \to b \leftarrow c$.

\begin{thmm}[{[{\citet{frydenberg1990}, Theorem~5.6}; {\citet{verma1990},
Theorem~1}]}]
\label{thmm:dag-equiv-immoral}
Directed acyclic graphs $\dG$ and $\dG'$ are Markov equivalent if
and only if they have the same skeleton and the same immoralities.
\end{thmm}

\subsection{Essential graphs}
\label{sec:essential-graphs}

An edge of a DAG $\dG$ is \emph{essential} if it has the same
direction in all Markov equivalent DAGs. The \emph{essential graph}
of $\dG$ is the graph in which all nonessential edges are replaced by
undirected edges.

Although not explored further in this work, the essential graph is a
type of \emph{chain graph}, a class of graphs that may have both
directed and undirected edges. For further details on chain graphs,
in particular their Markov properties and how they relate to
undirected and directed acyclic graphs, see \citet{frydenberg1990} and
\citet{andersson1997a}.

\begin{thmm}[{[{\citet{andersson1997}, Proposition~4.3}]}]
\label{thmm:dag-equiv-essential}
Directed acyclic graphs $\dG$ and $\dG'$ are Markov equivalent if
and only if they have the same essential
graph.
\end{thmm}

Unfortunately, there is no simple criterion for determining whether or
not an edge of a given DAG is essential, although
\citet{andersson1997} developed an iterative algorithm. This limits
their usefulness.

\subsection{Covered edge reversals}
\label{sec:cover-edge-revers}

A convenient characterisation of Markov equivalence can be given in
terms of edge reversals. An edge $a \to b$ of a DAG $\dG$ is
\emph{covered} if $\pa(b) = \pa(a) \cup\{a\}$.

\begin{thmm}[{[{\citet{chickering1995}, Theorem~2}]}]
\label{thmm:dag-equiv-reversal}
Directed acyclic graphs $\dG$ and $\dG'$ are Markov equivalent if
and only if there exists a sequence of DAGs
\[
\dG=\dG_0, \dG_1, \ldots, \dG_{k-1},
\dG_k= \dG'
\]
such that each $(\dG_{i-1},\dG_i)$ differ only by the reversal of
one covered edge.
\end{thmm}

This result is particularly useful for identifying properties that are
preserved under Markov equivalence, as it is only necessary to show
that the property is preserved under a covered edge reversal.

\subsection{Standard imset}
\label{sec:mark-equiv-imsets}

Imsets for undirected decomposable graphs were briefly mentioned in
Section~\ref{sec:clique-vector}. This formalism can be extended to
directed acyclic graphs. The \emph{standard imset} of a directed
acyclic graph $\dG$ is [\citet{studeny2005}, page 135]
\[
u_{\dG} = \delta(V) - \delta(\varnothing) + \sum
_{v \in V} \bigl[ \delta\bigl(\pa_{\dG} (v)\bigr) -
\delta\bigl(\pa_{\dG} (v) \cup\{v\} \bigr) \bigr],
\]
where $\delta(A) = (\ind_{S = A})_{S \subseteq V}$.

\begin{thmm}[{[\citet{studeny2005}, Corollary~7.1]}]
\label{thmm:dag-equiv-imset}
Directed acyclic graphs $\dG$ and $\dG'$ are Markov equivalent if
and only if $u_{\dG} = u_{\dG'}$.
\end{thmm}

\citet{studeny2009} give details of the relationship between the
imset and the essential graph of a DAG, and how one may be obtained
from the other.


\section{Proofs}
\label{sec:proofs}

\begin{pf*}{Proof of Theorem~\ref{thmm:sm-mark}}
The Markov property states that under $[\glaw, \famtheta]$,
%
\begin{equation}
\label{eq:sm-mark-1} X_A \indep X_B \bigmid
X_{A \cap B}, \rG, \bigl\{\rG\in\decom{A} {B}\bigr\}.
\end{equation}
Since if $\rG\in\decom{A}{B}$, then $\rG\simeq(\rG_A,\rG_B)$, we
can rewrite \eqref{eq:sm-mark-1} as
%
\begin{equation}
\label{eq:sm-mark-2} X_A \indep X_B \bigmid
X_{A \cap B}, \rG_A, \rG_B, \bigl\{\rG\in \decom{A} {B}
\bigr\}.
\end{equation}
As a consequence of Proposition~\ref{prp:compat-ci},
%
\begin{equation}
\label{eq:sm-mark-3} X_A \indep\rG_B \bigmid
X_{A \cap B} , \rG_A, \bigl\{\rG\in\decom {A} {B}\bigr\} ,
\end{equation}
and combined with \eqref{eq:sm-mark-2},
%
\begin{equation}
\label{eq:sm-mark-4} X_A \indep(X_B, \rG_B)
\bigmid X_{A \cap B}, \rG_A, \bigl \{\rG\in \decom{A} {B}\bigr\}.
\end{equation}
Furthermore, by the structural Markov property and Proposition~\ref{prp:compat-ci},
%
\begin{equation}
\label{eq:sm-mark-5} \rG_A \indep(X_B, \rG_B)
\bigmid\bigl \{\rG\in\decom{A} {B}\bigr\},
\end{equation}
and we can further condition on $X_{A \cap B}$. The result follows
from this and \eqref{eq:sm-mark-4}.
\end{pf*}

\begin{pf*}{Proof of Theorem~\ref{thmm:struct-markov-exponential}}
For any $C \subseteq V$, define $\G^{(C)}$ as in the proof of
Theorem~\ref{thmm:tform}, and let $\glaw$ have density $\pi$.

Suppose that $\glaw$ is structurally Markov. For any $\G\in
\sudg$, let $C_1,\ldots,C_k$ be a perfect ordering of the cliques,
and let $S_2, \ldots, S_k$ be the corresponding separators, and $H_i
= C_1 \cup\cdots\cup C_i$. Furthermore, recursively define the
graphs
\[
\G^{*(j)} = \cases{ \G^{(C_1)}, &\quad  $\mbox{if $j=1$}$, \vspace*{2pt}
\cr
\G^{*(j-1)}_{H_{j-1}} \gtimes\G^{(C_j)}_{(V \setminus H_{j-1})
\cup S_j}, &\quad  $
\mbox{if $j=2, \ldots, k$}$. }
\]
By Proposition~\ref{prp:sm-ratio}, for each $j=2, \ldots, k$,
\[
\pi \bigl( \G^{*(j)} \bigr) \pi \bigl( \G^{(S_j)} \bigr) = \pi
\bigl( \G^{*(j-1)} \bigr) \pi \bigl( \G^{(C_j)} \bigr).
\]
Note that $\G^{*(k)} = \G$. Then, by induction,
\[
\pi(\G) = \frac{
\prod_{j=1}^k \pi ( \G^{(C_j)}  )
}{
\prod_{j=2}^k \pi (\G^{(S_j)} )
} \propto\exp\bigl\{ \omega\cdot t(G) \bigr\}
\]
by Theorem~\ref{thmm:tform}, where $\omega_C = \log\pi
 (\G^{(C)} )$.

To show the converse let $(\omega)_A = (\omega_S)_{S \subseteq A}$.
By Lemma~\ref{lem:ctdecom},
\begin{eqnarray*}
&&\pi \bigl(\G_A \tmid\G_B, \bigl\{\G\in\decom{A} {B}\bigr\}
\bigr)
\\
&&\qquad \propto\exp \bigl\{ (\omega)_A \cdot t(\G_A) + (
\omega)_B \cdot t(\G_B) - (\omega)_{A \cap B} \cdot
t(\G_{A \cap B}) \bigr\}
\\
&&\qquad \propto\exp \bigl\{ (\omega)_A \cdot t(\G_A) - (
\omega)_{A \cap
B} \cdot t(\G_{A \cap B}) \bigr\}
\\
&&\qquad \propto\pi \bigl( \G_A \tmid\bigl\{\G\in\decom{A} {B}\bigr\} \bigr).
\end{eqnarray*}
\upqed\end{pf*}

\begin{pf*}{Proof of Theorem~\ref{thmm:dsmp-exp}}
If the law is in the exponential family \eqref{eq:dsmp-exp}, then,
by Theorem~\ref{thmm:d-clique-ancestral},
\[
\pi\bigl(\D| \ances(A)\bigr) \propto\exp\bigl\{ \omega\cdot\bigl[t(
\D_A) + t(\D _{V|A})\bigr] - \omega_A \bigr\}
\propto p\bigl(\D_A|\ances(A)\bigr) p\bigl(\D_{V|A} |
\ances(A)\bigr),
\]
and hence the law must be structurally Markov.

For the converse, define $\D^{(A)}$ to be the dagoid in which the
induced dagoid on $A \subseteq V$ is complete, but otherwise sparse
(in other words, the remainder dagoid, $\D^{(\varnothing)}_{V|A}$, of
the sparse dagoid $\D^{(\varnothing)}$ corresponding to complete
independence).

Select some $\dG\in\D$, and let $v_1,\ldots,v_d$ be a well
-ordering of $V$. Recursively define the dagoids
\[
\D^{*(i)} = \cases{ \D^{(\{v_1\})}, &\quad  $ \mbox{if $i=1$,}$
\vspace*{2pt}
\cr
\D^{*(i-1)}_{\pr(v_i)} \dprod\D^{(\{v_i\} \cup
\pa(v_i))}_{v_i|\pr(v_i)},
&\quad  $\mbox{otherwise.}$}
\]

By Proposition~\ref{prp:dag-sm-ratio}, for $i=2,\ldots,d$,
\[
\pi\bigl(\D^{*(i-1)}\bigr) \pi\bigl(\D^{(\{v_i\} \cup\pa(v_i))}\bigr) = \pi\bigl(
\D^{*(i)} \bigr) \pi\bigl( \D^{(\{v_i\} \cup\pa(v_i))}_{\pr(v_i)} \dprod
\D^{*(i-1)}_{v_i|\pr(v_i)} \bigr).
\]
However,
\[
\D^{(\{v_i\} \cup\pa(v_i))}_{\pr(v_i)} \dprod\D^{*(i-1)}_{v_i|\pr
(v_i)} =
\D^{(\pa(v_i))}.
\]
Therefore, since $\D^{*(d)} = \D$,
\[
\pi(\D) = \Biggl[ \prod_{i=1}^d \pi\bigl(
\D^{(\{v_i\} \cup\pa(v_i))}\bigr) \Biggr] \bigg/ \Biggl[ \prod_{i=2}^d
\pi\bigl(\D^{( \pa(v_i))}\bigr) \Biggr],
\]
which is of the form in \eqref{eq:dsmp-exp} with
\[
\omega_A = \log\pi\bigl( \D^{(A)} \bigr). 
\]
\upqed\end{pf*}
\end{appendix}



%





\printaddresses
\end{document}